\documentclass[final,3p]{elsarticle}

\usepackage{lineno}
\usepackage[colorlinks=true,breaklinks=true,pdftex]{hyperref}
\modulolinenumbers[1]

\usepackage[english]{babel}
\usepackage[utf8]{inputenc}
\usepackage{typearea}
\usepackage{latexsym}
\usepackage{lineno,hyperref}
\usepackage{amsmath, amssymb, mathtools,amsthm}
\usepackage{algorithm}
\usepackage{algorithmic}
\usepackage{graphicx}
\usepackage[abs]{overpic}
\usepackage{subcaption}
\usepackage{nicefrac}
\usepackage{textcomp} 
\usepackage[export]{adjustbox}

\newtheorem*{conjecture*}{Conjecture}
\newtheorem{remark}{Remark}

\DeclareRobustCommand{\ShowColormap}{\raisebox{-0.14em}{\includegraphics[height=.8em]{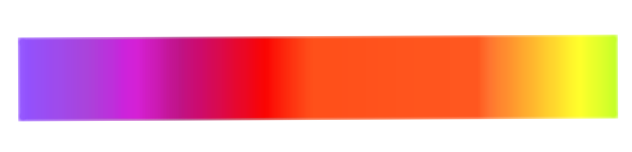}}}

\DeclareMathOperator*{\argmin}{arg\,min}
\newcommand{\rev}[1]{}

\newcommand{\dd}{\textnormal{d}}
\newcommand{\dist}{\textnormal{dist}}
\newcommand{\grad}{\textnormal{grad}}

\newcommand{\fcdot}{\; \cdot\; }

\newcommand*{\quark}{\setbox0\hbox{$x$}\hbox to\wd0{\hss$\cdot$\hss}}

\usepackage{xcolor}
\newcommand{\todo}[1]{}

\newtheorem{theorem}{Theorem}

\newtheorem{definition}{Definition}[section]
\newtheorem{corollary}{Corollary}[section]

\journal{CAGC}

\biboptions{authoryear}

\let\cite\citep

\begin{document}
	
	\begin{frontmatter}
		
		\title{
			De Casteljau's Algorithm in Geometric Data Analysis: Theory and Application
		}
		\author[first]{Martin Hanik\corref{cor}}
		\cortext[cor]{Corresponding author}
		\ead{hanik@zib.de}
		\author[last]{Esfandiar Nava-Yazdani}
		\ead{navayazdani@zib.de}
		\author[last]{Christoph von Tycowicz}
		\ead{vontycowicz@zib.de}
		\address[first]{Freie Universität Berlin, Germany}
		\address[last]{Zuse Institute Berlin, Germany}

		\begin{abstract}
			For decades, de Casteljau's algorithm has been used as a fundamental building block in curve and surface design and has found a wide range of applications in fields such as scientific computing, and discrete geometry to name but a few. With increasing interest in nonlinear data science, its constructive approach has been shown to provide a principled way to generalize parametric smooth curves to manifolds. These curves have found remarkable new applications in the analysis of parameter-dependent, geometric data. This article provides a \rev{survey} of the recent theoretical developments in this exciting area as well as its applications in fields such as geometric morphometrics and longitudinal data analysis in medicine, archaeology, and meteorology. 
			
		\end{abstract}
		
		
	\end{frontmatter}
	
	
	\section{Introduction}
	When what we call today B\'{e}zier curves were first studied in 1912 by Sergei Natanovich Bernstein~\cite{bernstein1912demo} as a means to function approximation, they were a purely theoretical tool. 
	This changed about 40 years later when Paul de Faget de Casteljau and Pierre Étienne Bézier\footnote{The curves bear the name of Pierre Étienne Bézier, as his employer, the car manufacturer Renault, allowed him to publish his works. Citroën, where Paul de Faget de Casteljau was employed, had a more restrictive policy, and his contributions were only made public several years after Bézier's.}, almost simultaneously but independently, searched for mathematical tools that would provide the means to intuitively construct and manipulate complex shapes using the still novel digital computers~\cite{farouki2012bernstein}.
	Both chose to use what are today called Bézier curves. With their elegant geometric properties, they meet the requirements; furthermore (as was found out later) they exhibit excellent numerical stability~\cite{farouki1987numerical}.
	The influence of the work of de Casteljau and Bézier was tremendous, and, today, Bézier curves are used in a multitude of applications. Especially in computer graphics and computer-aided geometric design, they are indispensable~\cite{farin2002curves}. 
	
	One of de Casteljau's great contributions is the algorithm that bears his name today. It is the fundamental tool to compute Bézier curves from a finite number of control points. (\citet{boehm1999casteljau} provide an introduction to the algorithm with historical remarks; a short autobiographical sketch by de Casteljau himself about the time leading to its invention can be found in~\cite{de1999casteljau}.) Given $k+1$ control points $p_0,\dots,p_k \in \mathbb{R}^d$, the value of the corresponding Bézier curve $\beta$ at $t \in [0,1]$ is computed iteratively according to the rule
	\begin{align}
		\label{def:euk_deCasteljau}
		\begin{split}
			&\beta_i^0(t):=p_i, \\
			&\beta_i^r(t):=(1-t)\beta_i^{r-1}(t) + t\beta_{i+1}^{r-1}(t), \quad r=1,\dots,k \quad i=0,\dots,k-r.
		\end{split}
	\end{align} 
	The resulting curve $\beta$ is a polynomial curve with coefficients $p_0, \dots, p_k$ when expanded in the Bernstein basis; i.e., with basis polynomials
	$b^k_j :[0,1] \to [0,1]$ defined by
	\begin{equation*}
		b^k_j(t) := \binom{k}{j}(1-t)^{k-j}t^j, \quad j=0,\dots,k,
	\end{equation*}
	it is of the form
	\begin{equation*}
		\beta(t) = \sum_{j=0}^k p_j b^k_j(t).
	\end{equation*}
	Although $\beta$ is written in terms of $k+1$ polynomials, its actual degree may be less than $k$~\cite[Sec.\ 5]{farouki2012bernstein}. 
	The algorithm can be extended to build triangular and quadrangular nets~\cite{de1963courbes,boehm1999casteljau}. Furthermore, one can string together Bézier curves to form Bézier splines, and simple geometric conditions determine the smoothness of the composite curve. 
	Another extension of de Casteljau's algorithm yields so-called polar forms or blossoms~\cite{de1985mathematiques,ramshaw1989blossoms}, which are particularly useful in the theoretical analysis of Bézier splines and algorithms for their computation. 
	
	Interest in modeling \textit{manifold}-valued curves (such as trajectories in the Lie group of rigid body motions) emerged roughly 30 years ago. It was then discovered that de Castejau's algorithm can be generalized from flat Euclidean space to curved manifolds~\cite{ParkRavani1995,crouch1999casteljau,lin2001cagd}.
	The necessary observation is that (\ref{def:euk_deCasteljau}) can be seen as following a straight line from $\beta_i^{r-1}$ to $\beta_{i+1}^{r-1}$ for a fixed amount of time $t$. As straight lines are the geodesics of Euclidean space, one can follow the geodesics of a curved space to generalize the algorithm.
	Although there is no analog to the Bernstein basis in manifolds (which facilitates theoretical analyses in the Euclidean case), it can be shown that the resulting \textit{generalized Bézier curves} inherit many geometric properties their Euclidean special cases possess~\cite{PN2007, NP2013}, including the easy construction of differentiable Bézier splines. One of the most important features of generalized Bézier curves is the (generalized) de Casteljau algorithm itself: Whenever, in a given manifold, geodesics can be computed efficiently, this is also true for generalized Bézier curves.     
	
	Several traditional applications of Bézier curves can also be found for their manifold-valued generalizations: They have been used, e.g., for curve modeling~\cite{morera2008modeling, SharpCrane2020, Macinelli_ea2023} and computer animation through interpolation~\cite{ParkRavani1995, GousenbourgerSamirAbsil2014, rumpfbezierimage, Brandt2016, SAMIR2019371}.
	
	
	To an increasing extent, manifold-valued Bézier curves have been used in \textit{data analytic} tasks---applications that de Casteljau probably did not envision when he derived his algorithm. Indeed, a wide range of applications in areas like  morphology~\cite{Kendall_ea2009, DrydenMardia2016}, action recognition~\cite{ParkBobrowPloen1995,  VeeriahZhuanQi2015, VemulapalliChellapa2016, Huang_ea2017}, and medical imaging~\cite{Arsigny_ea2007, BauerBruverisMichor2014, PennecSommerFletcher2020, AmbellanZachowTycowicz2021, Ambellan2022} work with manifold-valued data to the extent that relying on a (possibly hard to find) Euclidean embedding space does not suffice.
	Several data fitting approaches relying on generalized Bézier splines have been proposed and tested on different types of data~\cite{BergmanGousenbourger2018, GousenbourgerMassartAbsil2019, BakShinKoo2023}.
	
	A methodically related line of research has focused on the \textit{statistical} analysis of manifold-valued data using generalized Bézier splines as the underlying model. As one of the most fundamental statistical procedures, regression methods play an important role in the statistical analysis of the dependence of a response variable on an explanatory variable. 
	Parametric linear regression was first brought to manifolds in~\cite{Fletcher2013} by modeling trends as geodesics.
	The latter, while adequate in various circumstances, frequently fall short, e.g., when saturation effects or cyclic behavior are observed.
	Generalized Bézier splines, which are far more adaptable than geodesics, have successfully been utilized in~\cite{Hanik_ea2020} to capture such effects. 
	
	When the data has more complicated correlations, e.g., longitudinal data, hierarchical (or mixed-effects) models are needed to account for them. Extending the regression approach, a hierarchical model based on generalized Bézier splines was proposed that can capture a multitude of manifold-valued effects~\cite{HanikHegevonTycowicz2022}. An important concept underlying the model is that of \textit{Bézierfolds}: the manifolds consisting of all Bézier splines of the same type over a certain neighborhood. It was first discussed in~\cite{HanikHegevonTycowicz2022} and then further clarified in~\cite{NavaYazdaniAmbellanHaniketal2023}. Two different Riemannian metrics were proposed for Bézierfolds, namely, an ``integral'' metric~\cite{HanikHegevonTycowicz2022} that can be used for arbitrary splines and a generalization of the Sasaki metric~\cite{NavaYazdaniAmbellanHaniketal2023} for cubic splines.
	
	\rev{There are several closely related concepts for modeling smooth curves, such as B-splines and subdivision schemes.
		Replacing affine by geodesic averaging, subdivision schemes were generalized to manifolds~\cite{WALLNER2005593, NavaYu} and more recently applied to the modeling of curves in shape spaces~\cite{Brandt2016, HUBER2017313}.
		Additionally, spacetime contraints were used to compute optimal trajectories as generalizations of traditional cubic B-splines, named wiggly splines~\cite{kass2008animating}, and were adapted to curved shape spaces~\cite{schulz2014animating,schulz2015animating}.
		While these types of manifold-valued curves can exhibit similar flexible shapes, their computation requires expensive approximation methods.
		This also applies to curves that are derived from a variational setting on Riemannian manifolds, such as Riemannian cubics~\cite{Hinkle2014,NoakesRatiu2016,camarinha2022riemannian}.}
	
	The main goal of this article is to give an overview of \rev{manifold-valued Bézier splines in data analysis} and, in doing so, to provide an introduction to an intriguing subfield that builds upon de Casteljau's research. It is thus, above all, a \rev{survey} article. Nevertheless, several presented results have appeared only in one author's Ph.D.\ thesis~\cite{Hanik_phdthesis}.
	\rev{In particular, the notion of Bézierfolds (Riemannian manifolds of Bézier splines) and their essential properties (i.e., Thm.~\ref{thm:bezierfold} and Cor.~\ref{thm:bezierfold_spline}), have not been established before aside from the cubic case. Additionally, the discussion on the equivalence of maximum likelihood and least squares estimation for spline regression in manifolds is new; this includes Thm.~\ref{thm:max_like}.}
	
	As a disclaimer, we emphasize that the exposition in this paper should not be read as an objective and balanced overview of applications of generalized Bézier curves as a whole, but rather as a personal perspective biased by the authors' research interests.
	
	This article is structured as follows. Sec.~\ref{sec:splines} starts with a short summary of the necessary mathematical facts from differential geometry. Then, we define Bézier splines in Riemannian manifolds and collect some of their properties. In Sec.~\ref{bezierfolds}, we prove that there is a smooth manifold structure for certain sets of Bézier splines and discuss 2 possible choices of a Riemannian metric for them. We also state a conjecture that the results hold more generally. Sec.~\ref{sec:spline_regression} shows how Bézier splines can be used for highly flexible regression in manifolds. In Sec.~\ref{sec:hierarchical_model}, we then show how the results from Sec.~\ref{bezierfolds} are utilized for the hierarchical modeling of manifold-valued data. Last but not least, applications to real-world data are summarized in Sec.~\ref{sec:applications}. 
	
	It is not necessary to read all sections linearly. Sec.~\ref{sec:splines} is the only necessary prerequisite for the regression scheme presented in Sec.~\ref{sec:spline_regression}. The hierarchical model from Sec.~\ref{sec:hierarchical_model}, on the other hand, builds on all prior sections.
	
	\section{Manifold-valued Bézier Splines} \label{sec:splines}
	In this section, we recall manifold-valued Bézier splines and discuss their properties. In particular, we show that certain sets of Bézier curves constitute manifolds; two Riemannian metrics---one that is generally usable and one for the cubic case---are also given. We start with a brief summary of the necessary facts from Riemannian geometry; a good reference on this is~\cite {doCarmo1992}. We always use ``smooth'' for ``infinitely often differentiable.'' 
	
	\subsection{Background: Riemannian Geometry}
	A Riemannian manifold is a differentiable manifold $M$ with a Riemannian metric $\langle \fcdot, \fcdot \rangle$ that, for each $p \in M$, assigns a smoothly varying scalar product $\langle \fcdot, \fcdot \rangle_p$ to the tangent space $T_p M$. The metric induces a distance function $\dist$ on $M$.   
	If $f: M \to N$ is a smooth map between two manifolds $M$ and $N$, we denote the derivative of $f$ at $p$ in direction $v \in T_pM$ by $\textnormal{d}_pf(v)$. It yields a linear mapping $\textnormal{d}_pf: T_pM \to T_{f(p)}N$. 
	
	One of the fundamental objects every manifold possesses is the so-called Levi-Civita connection. Given two vector fields $X, Y$ on $M$, the latter allows to differentiate $Y$ along $X$, resulting in a new vector field $\nabla_X Y$. For the special case of Euclidean space, the Levi-Civita connection is simply the directional derivative. 
	
	The connection induces the Riemannian curvature tensor $R$, which intuitively measures the local deviation of $M$ from flat space by quantifying the failure of ``partial derivatives'' to commute.
	It takes three vector fields $X, Y, Z$ as input and is defined by
	$R(X, Y) Z := \nabla_X \nabla_Y Z - \nabla_Y \nabla_X Z - \nabla_{[X,Y]}Z,$
	where $[\fcdot, \fcdot]$ denotes the Lie bracket of vector fields; see \cite[Ch. 4]{doCarmo1992} for more details. 
	
	A vector field $X$ along $\alpha$ is called parallel if $\nabla_{\alpha'}X = 0;$
	we also say that $X_{\alpha(b)}$ is the parallel transport of $X_{\alpha(a)}$ to $T_{\alpha(b)}M$.
	
	A geodesic $\gamma$ is a generalized straight line. Denoting its derivative by $\gamma':= \frac{\dd}{\dd t} \gamma$, it is a curve without acceleration, i.e., $\nabla_{\gamma'}\gamma' = 0$. In Euclidean space, geodesics are simply straight lines that are traced with constant velocity. 
	
	It will be crucial for us that every point in $M$ has a normal convex neighborhood\footnote{Other names are also used for such a neighborhood; (at least) ``geodesically convex'' and ``strongly convex'' neighborhood can also be found.} $U$, in which each pair $p,q \in U$ is joined by a unique length-minimizing geodesic $[0,1] \ni t \mapsto \gamma(t;p,q)$ that never leaves $U$.
	Each geodesic is also differentiable w.r.t.\ its start and endpoints. 
	We will fix some normal convex neighborhood $U$ in the following.
	
	
	The Riemannian exponential is another very important map. Let $v \in T_pM$ and $\gamma$ be a geodesic in $U$ with $v = \gamma'(0;p,q)$. Then, $\exp_p(v) := q$ defines the exponential map at $p$. It has an inverse $\log_p$ called the logarithm at $p$, which fulfils $\log_p(q) = \gamma'(0;p,q)$. 
	Thus, the derivative $(\dd\exp_{p})_{v}: T_vT_pM \cong T_pM \to T_{\exp_p(v)}M$ of $\exp_p$ at $v$ is invertible for all $v \in \log_p(U)$. 
	
	Using the metric we can define the gradient of a smooth function $f: M \to \mathbb{R}$ implicitly by requiring
	\begin{equation} \label{eq:gradient}
		\dd f_p(v) = \langle \textnormal{grad}_pf, v \rangle_p,
	\end{equation}
	for all $v \in T_pM$ and all $p \in M$, generalizing the ordinary Euclidean gradient. Importantly, just like in the Euclidean case, $\dd_p f = 0$ if and only if $\grad_p f = 0$. Furthermore, the gradient still points in the direction of the steepest ascent, so optimization algorithms relying on gradients can be transferred to Riemannian manifolds~\cite{AbsilMahonySepulchre2007}.
	
	
	Finally, the mean as the most fundamental statistic of a data set can be generalized to manifolds in the form of the (sample) Fréchet mean. Given data $q_1,\dots,q_n \in U$, it is defined by
	\begin{equation} \label{def:frechet_mean}
		\overline{q} := \argmin_{p \in U}\sum_{i=1}^n \dist(q, q_i)^2.
	\end{equation}

	\subsection{Bézier Splines}
\begin{figure}[t]
		\centering
		\includegraphics[width=.6\linewidth]{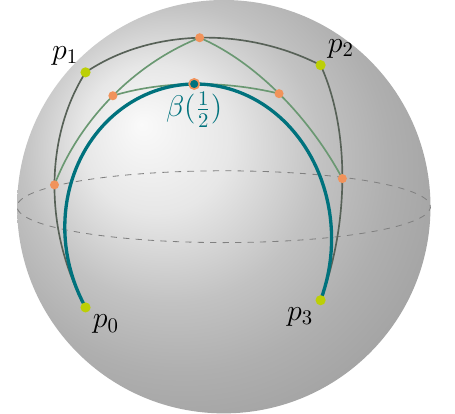}
		\caption{Cubic Bézier curve $\beta$ on the sphere $\mathcal{S}^2$ and the construction of $\beta(\nicefrac{1}{2})$ by the de Casteljau algorithm.}
		\label{fig: Bezier_curve_sphere}
	\end{figure}

	We can now define Bézier curves in $M$ by generalizing de Casteljau's algorithm~(\ref{def:euk_deCasteljau}).
	\begin{definition}
		[De Casteljau's Algorithm on Manifolds]\label{dc}
		Let control points $p_0,\dots,p_k \in U$ be given. For all $t \in [0,1]$, we set
		\begin{align*}
			&\beta_i^0(t):=p_i, \\
			&\beta_i^r(t):=\gamma(t; \beta_i^{r-1}(t), \beta_{i+1}^{r-1}(t)), \quad r=1,\dots,k, \quad  i=0,\dots,k-r.
		\end{align*}
		We call $\beta := \beta^k_0 : [0,1] \to U$ \textnormal{Bézier curve of degree $k$ with control points} $p_0,\dots,p_k$.  
	\end{definition}
	We also write $\beta(t; p_0,\dots,p_k)$ instead of $\beta(t)$ whenever we want to make the dependence on the control points explicit. The algorithm is visualized for the 2-dimensional sphere $\mathcal{S}^2$ in Fig.~\ref{fig: Bezier_curve_sphere}.
	
	The following theorem collects important properties of generalized Bézier curves from~\cite{PN2007}.
	\begin{theorem} \label{thm: properties_bezier}
		Let $\beta$ be a Bézier curve of degree $k$ with control points $p_0,\dots,p_k \in U$. Then,
		\begin{itemize}
			\item [(i)] $\beta(0) = p_0$ and $\beta(1) = p_k$;
			\item[(ii)] $\beta'(0) = k \log_{p_0}(p_1)$ and $\beta'(1) = -k \log_{p_k}(p_{k-1})$;
			\item[(iii)] $\nabla_{\beta'} \beta'(0) = k(k-1) \widetilde{u}_0$, where
			$$\widetilde{u}_0 := \begin{cases} \gamma'(0; p_1, p_2), & \textnormal{if } p_0=p_1, \\
				(\dd \exp_{p_0})^{-1}_{\gamma'(0;p_0,p_1)} \big(\gamma'(0; p_1, p_2) - \gamma'(1;p_0,p_1) \big), & \textnormal{if } p_0 \ne p_1;\end{cases}$$
			\item[(iv)] $\nabla_{\beta'} \beta'(1) = k(k-1) \widetilde{u}_k$, where
			$$\widetilde{u}_k := \begin{cases} -\gamma'(1; p_{k-2}, p_{k-1}), & \textnormal{if } p_{k-1}=p_k, \\
				(\dd \exp_{p_k})^{-1}_{-\gamma'(1;p_{k-1},p_k)} \big(\gamma'(0; p_{k-1}, p_k) - \gamma'(1;p_{k-2},p_{k-1}) \big), & \textnormal{if } p_{k-1} \ne p_k.\end{cases}$$
		\end{itemize}
	\end{theorem}
	
	The theorem yields the following corollary.
	\begin{corollary}\label{cor: fix_control_points}
		Let $\beta$ be a Bézier curve of degree $k$ with control points $p_0,\dots,p_k \in U$. Then,
		\begin{itemize}
			\item[(i)] $p_1 = \exp_{p_0} \left(\frac{1}{k}  \beta'(0) \right)$ and $p_{k-1} = \exp_{p_k}\left(-\frac{1}{k}  \beta'(1) \right)$;
			\item[(ii)] $p_2 = \exp_{p_1} \left(\frac{1}{k(k-1)} w_0\right)$
			and $p_{k-2} = \exp_{p_{k-1}} \left( \frac{1}{k(k-1)} w_k\right)$,
			where $u_0:= \nabla_{\beta'}\beta'(0)$, $u_k := \nabla_{\beta'}\beta'(1)$ and
			\begin{align*}
				w_0 &:= \begin{cases}
					u_0, & \textnormal{if } p_0 = p_1, \\
					(\dd \exp_{p_0})_{\gamma'(0;p_0,p_1)} \big(u_0 + k(k-1) \gamma'(0; p_0, p_1) \big), & \textnormal{if } p_0 \ne p_1,
				\end{cases} \\
				w_k &:= \begin{cases}
					u_k, & \textnormal{if } p_{k-1} = p_k, \\
					(\dd \exp_{p_k})_{-\gamma'(1;p_{k-1},p_k)} \big(u_k - k(k-1) \gamma'(1; p_{k-1}, p_k) \big), & \textnormal{if } p_{k-1} \ne p_k.
				\end{cases}
			\end{align*}
		\end{itemize}
	\end{corollary}
	
	These properties can be exploited to define $C^1$ or $C^2$ B\'ezier splines~\cite{PN2007, gousenbourger2020interpolation}. We give the definition for the $C^1$ case.
	\begin{definition} \label{def: spline}
		For $i=0,\dots,L-1$, let $p^{(i)}_0, \dots, p^{(i)}_{k_i} \in U$ be the control points of $L\ge2$ cubic Bézier curves $\beta^{(0)},\dots,\beta^{(L-1)}$ such that 
		\begin{equation} \label{eq:C1_condition}
			p^{(i)}_{k_i} = p^{(i+1)}_0 \quad \text{and} \quad  \gamma \left( \frac{k_{i+1}}{k_i + k_{i+1}}; p^{(i)}_{k_i -1}, p^{(i+1)}_1 \right) = p^{(i)}_{k_i}
		\end{equation}
		for all $i=0,\dots,L-2$.
		The \textnormal{($k_1,\dots,k_L$)-Bézier spline $B$ with control points} 
		$$\left( p^{(i)}_0,\dots,p^{(i)}_{k_i} \right)_{i=0,\dots,L-1}$$ 
		is then defined by
		\begin{equation*}
			B(t) := \begin{cases} \beta^{(0)}\left(t;p^{(0)}_0,p^{(0)}_1,p^{(0)}_2,p^{(0)}_3\right), &\quad t \in [0,1],\\ \beta^{(i)} \left(t-i;p^{(i)}_0,p^{(i)}_1,p^{(i)}_2,p^{(i)}_3 \right), &\quad t \in (i,i+1], \quad i=1,\dots,L-1. \end{cases}
		\end{equation*}
	\end{definition}
	When clear from context, we do not mention the degrees.
	It follows directly from Thm.~\ref{thm: properties_bezier} and Cor.~\ref{cor: fix_control_points} that Bézier splines are $C^1$.
	We also consider Bézier curves as splines with a single segment (the $L=1$ case).
	
	We can also ensure that a Bézier spline is cyclic. Indeed a Bézier spline is closed if and only if (\ref{eq:C1_condition}) extends cyclically; this holds if additionally
	\begin{equation} \label{eq:C1_condition_cyclic}
		p^{(L-1)}_k = p^{(0)}_0 \quad \text{and} \quad  \gamma \left(
		\frac{k_0}{k_0 + k_{L-1}};p^{(L-1)}_{k_{L-1}-1},p^{(0)}_1 \right) = p^{(0)}_0,
	\end{equation}
	Closed $C^1$ splines can readily be extended to $C^1$ cyclic curves defined on the whole of $\mathbb{R}$ when the input parameter is viewed modulo $L$. There are many applications where (quasi)-periodic behavior in manifolds is of central interest: E.g., the motion of a beating heart and gait patterns, which have been successfully modeled in curved shape space and the group of rigid-body motions, respectively.
	
	The fact that one can model periodic behavior with Bézier splines by simple conditions on their control points is a striking feature. Other generalizations of polynomial curves like Riemannian cubics~\cite{NoakesRatiu2016}, which are defined as solutions to initial value problems, do not allow for such easy modeling of cyclic trends.  
	
	Because of the restrictions (\ref{eq:C1_condition}) and (\ref{eq:C1_condition_cyclic}), not all control points are free parameters of a Bézier spline. Indeed, if $p^{(i)}_{k_i}=p^{(i+1)}_0$ are \textit{connecting control points}, joining the $i$-th to the $(i+1)$-th segment\footnote{For closed splines the segment numbers are meant modulo $L$, so the connection between the first and last segment is included.}, we omit $p^{(i+1)}_0$ as a variable. Furthermore, we \textit{choose} that $p^{(i)}_{k_i}$ and its \textit{predecessor} $p^{(i)}_{k_i-1}$ are free variables, while the \textit{successor} $p^{(i+1)}_1$ shall also be eliminated as a variable. With $\overline{k} :=(k_{i-1}+k_{i})/k_i$ we, therefore, have the additional condition $p^{(i+1)}_1 = \gamma(\overline{k};p^{(i)}_{k_i-1},p^{(i)}_{k_i})$ at the connection.
	Setting 
	\begin{equation} \label{def:K}
		\begin{cases} \mathring{K} := k_0+k_1+\cdots+k_{L-2}+k_{L-1} - L - 1, & \textnormal{if $B$ is closed,} \\
			K := k_0+k_1+\cdots+k_{L-2}+k_{L-1} - L + 1, & \textnormal{else,}\end{cases}
	\end{equation}
	it then follows that the number of \textit{independent} control points of a general $C^1$ Bézier spline $B$ is $K+1$, while for closed Bézier splines have $\mathring{K} + 1$.
	\begin{remark}
		Clearly, closed Bézier splines are Bézier splines as well. Therefore, when considering \textnormal{all} Bézier splines of given degrees, we \textit{include} closed splines and think of all $K +1$ control points as free parameters; only when we \textnormal{restrict} to closed splines, this number drops to $\mathring{K} + 1$.
	\end{remark}
	\begin{figure}[t]
		\begin{subfigure}[c]{0.49\textwidth}
			\includegraphics[width=.99\linewidth]{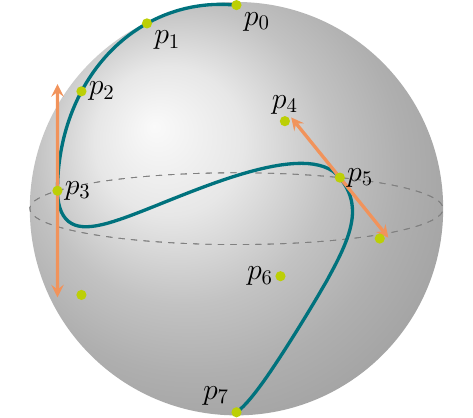}
			\subcaption{Bézier spline with 3 cubic segments on $\mathcal{S}^2$.}
			\label{fig:Bezier_spline_sphere}
		\end{subfigure}
		\hspace{2ex}
		\begin{subfigure}[c]{0.49\textwidth}
			\includegraphics[width=.99\linewidth]{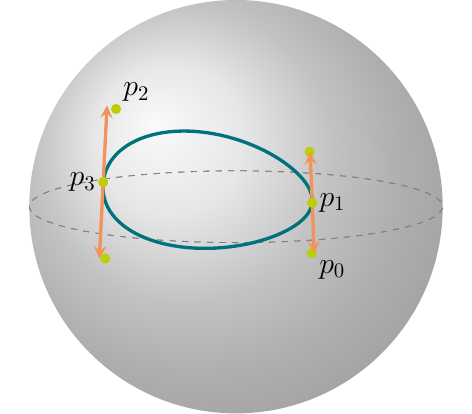}
			\subcaption{Closed spline with two cubic segments on $\mathcal{S}^2$.}
			\label{fig:closed_spline}
		\end{subfigure}
		\caption{Examples of Bézier splines. The unlabeled green points are dependent control points.}
	\end{figure}
	
	We denote the set of $K+1$ independent control points of a general spline $B$ by $p_0,\dots,p_K$, i.e.,
	$$(p_0,\dots,p_K) := \left(p^{(0)}_0,\dots,p^{(0)}_{k_0},p^{(1)}_2,\dots,p^{(1)}_{k_1},\dots,p^{(L-1)}_2,\dots,p^{(L-1)}_{k_{L-1}} \right) \in U^{K+1}.$$ 
	When restricting to closed splines, $p_0^{(0)}$ and $p_1^{(0)}$ are left out.
	An example of a non-closed $C^1$ spline with three cubic segments and nine independent control points is shown in Fig.~\ref{fig:Bezier_spline_sphere}, while a closed $C^1$ spline with two cubic segments and four independent control points is depicted in Fig.~\ref{fig:closed_spline}. 
	When it is important, we make the dependence of $B$ on its control points clear by writing $B(t;p_0,\dots,p_K)$ or $B(t;p_0,\dots,p_{\mathring{K}})$.

	\section{Spaces of B\' ezier Splines} \label{bezierfolds}
	We will now see that certain sets of Bézier splines over manifolds form manifolds themselves. To the best of our knowledge, this was first discussed in~\cite{HanikHegevonTycowicz2022}, albeit without rigorous proof; the latter was then given for cubic splines in~\cite{NavaYazdaniAmbellanHaniketal2023}. Here, we will extend the proof to splines with segments of maximal degree 5. After that, two Riemannian metrics for the resulting manifolds are discussed. The following background is the prerequisite for the Sasakian metric.
	
	\subsection{Background: Geometry of the Tangent Bundle} \label{sec:sasaki_metric}
	Next comes a summary of facts about the geometry of the tangent bundle; a good reference on the differential geometric properties of the tangent bundle is~\cite{Jost2017}; see~\cite{gudmundsson2002geometry} for a reference on the Sasaki metric.
	
	Let $M$ be a $d$-dimensional Riemannian manifold. The disjoint union $TM:= \bigsqcup_{p \in M} T_pM$ is the tangent bundle of $M$. It comes  with the bundle projection $\pi: TM \to M$, $(p,v) \mapsto p$.  
	The tangent bundle is a $2d$-dimensional manifold itself, and its natural Riemannian metric is the Sasaki metric. 
	To define the latter, we use the fact that the tangent bundle $TTM$ of $TM$ is the direct sum of a vertical subbundle $VTM$ (the kernel of the derivative $\textnormal{d}\pi$ of $\pi$) and a horizontal subbundle $HTM$ that is determined by the Levi-Civita connection of $M$ as vectors tangent to parallel vector fields; both $VTM$ and $HTM$ have rank $d$. Intuitively, horizontal vectors are directions in which only the footpoint changes, whereas the former is constant in vertical directions.  
	The Sasaki metric~\cite{sasaki1962differential} is the unique Riemannian metric $\langle \fcdot, \fcdot \rangle^S$ on $TM$ with the following properties:
	\begin{itemize}
		\item[(i)] The bundle projection $\pi$ is a Riemannian submersion, i.e.\ $\pi$ has maximal rank and $\textnormal{d}\pi$ preserves lengths of horizontal vectors.
		\item[(ii)] For any $p \in M$, the restriction of $\langle \fcdot, \fcdot \rangle^S$ to the tangent space $T_pM \subset TM$ coincides with $\langle \fcdot, \fcdot \rangle_p$.
		\item[(iii)] Let $u$ be a \textit{parallel} vector field along a curve $q: (-\varepsilon,\varepsilon) \to M$. Define $\zeta = (-\varepsilon,\varepsilon) \to TM$, $t \mapsto (q(t), u(t))$. Let further $p=q(0)$ lie on $q$ and $v: (-\varepsilon,\varepsilon) \to T_pM$ such that $u(0) = v(0)$. Define $\eta: (-\varepsilon,\varepsilon) \to TM$, $t \mapsto (p, v(t))$. Then $\dot{\zeta}(0) := \frac{d}{dt}\zeta(0) \in H_{(p,u(0))}TM$ and $\dot{\eta}(0) = \frac{d}{dt}\eta(0) \in V_{(p,u(0))}TM$ are orthogonal.
	\end{itemize}
	It follows that horizontal and vertical vectors are orthogonal.
	
	Let $(p,u) \in TM$. Because $H_{(p,u)}TM$ and $V_{(p,u)}TM$ are $d$-dimensional vector spaces, both can be identified with $T_pM$; in other words, we can view an element of $T_{(p,u)}TM$ as a tuple $(v,w) \in (T_pM)^2$. As horizontal and vertical components are orthogonal, the Sasaki metric between $(v_1,w_1), (v_2,w_2) \in T_{(p,u)}TM$ then reads
	\begin{equation} \label{def:Sasaki_metric}
		\langle (v_1,w_1), (v_2,w_2) \rangle^S_{(p,u)} = \langle v_1, v_2 \rangle_p + \langle w_1, w_2 \rangle_p.
	\end{equation}
	
	Geodesics of the Sasaki metric can be characterized in terms of geometric features of the underlying space.
	Let $\eta=(q,u): (-\varepsilon,\varepsilon) \to TM$ be a curve in $TM$ and $\dot{\eta} = (v,w)$. Denoting the Riemannian curvature tensor of $M$ by $R$, $\eta$ is a geodesic if (and only if) the coupled system
	\begin{align}
		\nabla_vv &= -R(u, w)v, \\
		\nabla_vw &= 0
		\label{eq:sasaki_geodesic}
	\end{align}
	holds. Intuitively, the above equations say that the footpoint curve $q$ bends according to the curvature of $M$, while the vector component $u$ changes at a constant rate.
	
	Further properties of $TM$ with the Sasaki metric $\langle \fcdot, \fcdot \rangle^S$, e.g.\ its curvature tensor and Levi-Civita connection, can be found in~\cite{gudmundsson2002geometry}.
	
	\subsection{Bézierfolds}
	We start with the definition of Bézierfolds.
	\begin{definition}[Bézierfolds]\label{def:bezierfold}
		Let $U$ be a normal convex neighborhood of a Riemannian manifold $M$. We define the \textnormal{Bézierfold of $(k_0,\dots,k_{L-1})$-Bézier splines over $U$} by
		\begin{align} \label{def:bezierfold}
			\begin{split}
				\mathcal{B}^{L}_{k_0,\dots,k_{L-1}}(U) := \{B: &[0,L]\to U\ \big|\ B \text{ is a $(k_0,\dots,k_{L-1})$-Bézier spline}\}.
			\end{split}
		\end{align}
		Its subset of closed splines  
		\begin{align} \label{def:closed_bezierfold}
			\begin{split}
				\mathring{\mathcal{B}}^L_{k_0,\dots,k_{L-1}}(U) := \{B: &[0,L]\to U\ \big|\ B \text{ is a closed $(k_0,\dots,k_{L-1})$-Bézier spline}\}.
			\end{split}
		\end{align}
		is called \textnormal{Bezierfold of closed $(k_0,\dots,k_{L-1})$-Bézier splines over $U$.}
	\end{definition}
	We often only say ``Bézierfold'' and write $\mathcal{B}(U)$ when the exact values of the parameters $(k_0,\dots,k_{L-1})$ are unimportant or clear.
	
	The following theorem shows that many Bézierfolds are not only arbitrary sets of functions but can be given the structure of finite-dimensional smooth manifolds.
	\begin{theorem} \label{thm:bezierfold}
		Let $M$ be a $d$-dimensional Riemannian manifold and $U \subseteq M$ a normal convex neighborhood. Furthermore, let $\mathcal{B}^1_k(U)$ and  $\mathring{\mathcal{B}}^1_k(U)$ be defined by (\ref{def:bezierfold}) and (\ref{def:closed_bezierfold}), respectively. Then, $\mathcal{B}^1_k(U)$ can be given the structure of a $(k+1)d$-dimensional smooth manifold for $k \in \{0,1,2,3,4,5\}$. Furthermore, with this structure, $\mathring{\mathcal{B}}^1_k(U)$ is a smooth $(k-1)d$-dimensional embedded submanifold of $B^1_k(U)$ for $k \in \{2,3,4,5\}$.
	\end{theorem}
	\begin{proof}
		Let $1 \le k \le 5$. We define the map
		\begin{align} \label{def:F_prime}
			F&:U^{k+1} \to \mathcal{B}^1_k(U), \\
			&(p_0,\dots,p_k) \mapsto \beta(\fcdot;p_0,\dots,p_k).
		\end{align}
		Thm.~\ref{thm: properties_bezier} and Cor.~\ref{cor: fix_control_points} imply that $F$ is injective and, thus, bijective. We can then push the product-manifold structure of $U^{k+1}$ forward along $F$ to $B^1_k(U)$ to turn it into a smooth manifold, thereby securing that $F$ is a smooth diffeomorphism~\cite[Ch.\ 30 § 9]{Postnikov2013}.\footnote{More precisely, we obtain a topology on $\mathcal{B}^1_k(U)$ by requiring that a subset is open if and only if its preimage under $F$ is open. A maximal atlas is then given analogously: If $(V,\phi)$ is a chart of $U^{k+1}$, then $(F(U), \phi \circ F^{-1})$ is a chart of $B^1_k(U)$ (and all charts are constructed this way).} Analogously, we can push the structure of $U^{k-1}$ to $\mathring{\mathcal{B}}^1_k(U)$ turning it into an embedded submanifold.
	\end{proof}
	This structure allows us to identify Bézier splines of low degrees and their control points.
	Application-wise, the spaces covered in Thm.~\ref{thm:bezierfold} are the most important ones as, traditionally, the degree of the underlying polynomial curve in regression is kept as low as possible to avoid overfitting, and this should not be different for regression in manifolds.
	Also, when curves of low degrees are not adequate, it can be advisable to use splines (built from them) instead of curves of higher degrees~\cite{MarshCormier2001}.
	
	From a theoretical perspective, it is nevertheless interesting to advance the result to higher degrees.
	Indeed, it seems very likely that $\mathcal{B}^1_k(U)$ and $\mathring{\mathcal{B}}^1_k(U)$ can always be given a manifold structure in the above way. We, therefore, state the following conjecture. 
	\begin{conjecture*}
		Thm.~\ref{thm:bezierfold} also holds for all $k \ge 6$.
	\end{conjecture*}
	We are unaware of explicit formulas for further control points, which would be necessary to extend our proof. 
	What is known is that the $\ell$-th covariant derivative of a Bézier curve at $t=0$ and $t=1$ is determined by the first and last $\ell+1$ control points, respectively~\cite[p.\ 112]{PN2007}; but not that a subset of them never suffices. 
	
	Of course, one could also try to use the inverse function theorem to show that $F$ is \textit{locally} invertible. Then, by restricting to a subset $\widetilde{U} \subseteq U$, we would obtain that $\mathcal{B}^1_k(\widetilde{U})$ is a manifold. To this end, delicate investigations of sums of Jacobi fields along different curves are necessary. We leave this for future work.
	
	Turning to Bézier splines, we can extend Thm.~\ref{thm:bezierfold}.
	\begin{corollary} \label{thm:bezierfold_spline}
		Let $M$ be a $d$-dimensional Riemannian manifold and $U \subseteq M$ a normal convex neighborhood. Furthermore, let $\mathcal{B}^L_{k_0,\dots,k_{L-1}}(U)$ and  $\mathring{\mathcal{B}}^L_{k_0,\dots,k_{L-1}}(U)$ be defined by (\ref{def:bezierfold}) and (\ref{def:closed_bezierfold}), respectively. Then, $\mathcal{B}^L_{k_0,\dots,k_{L-1}}(U)$ can be given the structure of a smooth manifold of dimension 
		$$\textnormal{dim} \left(\mathcal{B}^L_{k_0,\dots,k_{L-1}}(U) \right) = d \left( \sum_{i=0}^{L-1}k_i - L + 2 \right)$$
		if $k_i \in \{1,2,3,4,5\}$ for all $i=0,\dots, L-1$.
		Furthermore, if $k_i \in \{1,2,3,4,5\}$ for all $i=1,\dots,L-2$ such that $\sum_{i=0}^{L-1} k_i > L$, then $\mathring{\mathcal{B}}^L_{k_0,\dots,k_{L-1}}(U)$ is a smooth embedded submanifold in this structure. The dimension of $\mathring{\mathcal{B}}^L_{k_0,\dots,k_{L-1}}(U)$ is
		$$\textnormal{dim} \left(\mathring{\mathcal{B}}^L_{k_0,\dots,k_{L-1}}(U) \right) = d \left( \sum_{i=0}^{L-1} k_i - L \right).$$
	\end{corollary}
	\begin{proof}
		We extend $F$ so that
		\begin{align} \label{eq:diffeo_F}
			\begin{split}
				F&: U^{K+1} \to \mathcal{B}_{k_0,\dots,k_{L-1}}^L(U), \\ 
				&(p_0,\dots,p_K) \mapsto B(\fcdot;p_0,\dots,p_K), 
			\end{split}
		\end{align}
		where $K$ is given by Def.~\ref{def:K}.
		The assertion now follows using the arguments from the proof of Thm.~\ref{thm:bezierfold} for every segment. (Thereby, $F$ is also turned into a diffeomorphism.) 
		The construction for $\mathring{\mathcal{B}}^L_{k_0,\dots,k_{L-1}}(U)$ works analogously.
	\end{proof}
	
	\begin{remark} \label{rem:constant_bezierfolds}
		Several Bézierfolds coincide with sets of constant maps into $U$. Whereas for spaces that contain non-closed curves, this is (only) true for $\mathcal{B}^1_0(U)$, there are infinitely many spaces of closed curves, including $\mathring{\mathcal{B}}^1_0(U)$, $\mathring{\mathcal{B}}^1_1(U), \mathring{\mathcal{B}}^1_2(U), \mathring{\mathcal{B}}^2_{1,2}(U)$, $\mathring{\mathcal{B}}^3_{1,1,2}(U)$, etc. It is then clear that these spaces are all diffeomorphic to $U$.
	\end{remark}
	
	Since derivatives of diffeomorphisms map tangent spaces bijectively into each other, we can characterize the tangent spaces of a Bézierfold $\mathcal{B}_{k_0,\dots,k_{L-1}}^L(U)$ with the help of the map (\ref{eq:diffeo_F}): For each $B \in \mathcal{B}_{k_0,\dots,k_{L-1}}^L(U)$, we find
	\begin{align} \label{eq:tangent_space_primal}
		T_B\mathcal{B}_{k_0,\dots,k_{L-1}}^L(U) = \{ X&: [0,L] \to TM\ \big|\  \exists \ v_0 \in T_{p_0}M, \dots, v_K \in T_{p_K}M\  \forall\  t \in [0,L]: \nonumber\\
		&X(t) = \sum_{j=0}^K \dd_{p_j}B(t;p_0,\dots,p_{j-1},\fcdot,p_{j+1},\dots,p_K)(v_j)\}.
	\end{align}
	Every element $X \in T_B\mathcal{B}_{k_0,\dots,k_{L-1}}^L(U)$ is thus a vector field along $B$. We say that it is induced by the vectors $v_0, \dots, v_K$. The tangent spaces $T_B\mathring{\mathcal{B}}^L_{k_0,\dots,k_{L-1}}(U)$ of spaces with closed curves are given by the right-hand side of~(\ref{eq:tangent_space_primal}) when $K$ is replaced by $\mathring{K}$.
	Vector fields as in (\ref{eq:tangent_space_primal}) were studied by \citet{BergmanGousenbourger2018}. They derived explicit formulas for their computation if $M$ is a symmetric space, i.e., Riemannian manifolds for which reflections about points are isometries.
	
	In~\cite{NavaYazdaniAmbellanHaniketal2023}, a closely connected identification was used to give the space of \textit{cubic}\footnote{It is relatively straightforward to extend the construction to odd-degrees-only Bézierfolds, i.e.,  $\mathcal{B}^L_{k_0,\dots,k_{L-1}}(U)$ with $k_0,\dots,k_{L-1} \in \{1,3,5\}$.} Bézier splines a smooth manifold structure: The authors use the fact that an ordered pair $(p, q) \in U$ can be identified with $(p, \log_p(q)) \in TM$. A cubic spline $B \in \mathcal{B}^L_{3,\dots,3}(U)$ can then be encoded using the mapping
	\begin{align} \label{def:F}
		\begin{split}
			\widetilde{F}:\mathcal{B}^L_{3,\dots,3}(U) &\to (TU)^{L+1}, \\
			B &\mapsto \bigg( \big(B(0), \dot{B}(0) /3 \big), \dots, \big( B(L), \dot{B}(L) / 3 \big) \bigg);
		\end{split}
	\end{align}
	Thm.~\ref{thm: properties_bezier} and Cor.~\ref{cor: fix_control_points} yield
	\begin{align*}
		\widetilde{F}(B) = \Bigg( \bigg(p^{(0)}_0, \log_{p^{(0)}_0}\big(p^{(0)}_1\big) \bigg), \dots, \bigg(p^{(L-1)}_0, \log_{p^{(L-1)}_0}&\big(p^{(L-1)}_1\big) \bigg), \\
		&\bigg(p^{(L-1)}_3, -\log_{p^{(L-1)}_3}\big(p^{(L-1)}_2\big) \bigg) \Bigg).
	\end{align*}
	For splines in $\mathring{\mathcal{B}}^L_{3,\dots,3}(U)$, the last element is left out.
	
	The equivalence of the induced structures follows from the fact that the differentiable structure of the tangent bundle $TM$ is induced by the one of $M$. The mapping $\widetilde{F}$ is helpful in obtaining a Riemannian metric for cubic Bézierfolds as we will see shortly.
	
	The map $\widetilde{F}^{-1}$ can be used for another characterization of the tangent spaces (\ref{eq:tangent_space_primal}) for cubic Bézierfolds. We obtain
	\begin{align} \label{eq:tangent_space}
		\begin{split}
			T_B\mathcal{B}^L_{3,\dots,3}(U) = \{ X: [0,L] &\to TM\ \big|\  \exists \ (v_0, w_0) \in T_{(p_0, u_0)}TM, \dots, \\ 
			&(v_L, w_L) \in T_{(p_L, u_L)}TM: X = \sum_{i=0}^L \textnormal{d}_{(p_i,u_i)}\widetilde{F}^{-1}\big((v_i,w_i)\big)\},
		\end{split}
	\end{align}
	while for $T_B\mathring{\mathcal{B}}^L_{3,\dots,3}(U)$ the last summand of the vector fields is left out.
	The derivatives of $\widetilde{F}^{-1}$ and $F$ can be computed using the same tools (so-called Jacobi fields). In particular, there are explicit formulas available in symmetric spaces to compute the terms $\textnormal{d}_{(p_i,u_i)}\widetilde{F}^{-1}\big((v_i,w_i)$ in~(\ref{eq:tangent_space}); see~\cite{Fletcher2013, NavaYazdaniAmbellanHaniketal2023}.

	\subsection{Riemannian Metrics on Bézierfolds}
	
	As a manifold, a Bézierfold can be endowed with a Riemannian metric. In the following, we describe two different possible choices---one general and one for the special case of cubic Bézier splines. We state everything for spaces $\mathcal{B}(U)$ of general splines, but, clearly, the metrics can be restricted to submanifolds $\mathring{\mathcal{B}}(U)$ of closed curves.
	\subsubsection{Integral-based Metric}\label{sec:integral-based_metric}
	In~\cite{HanikHegevonTycowicz2022}, the authors proposed to endow $\mathcal{B}^L_{k_0,\dots,k_{L-1}}(U)$ with the metric of Srivastava and Klassen from~\cite[Sec.~3.3]{SrivastavaKlassen2016}. 
	\begin{definition} \label{eq:metric_bezierfold}
		Let $\langle \fcdot, \fcdot \rangle$ be the Riemannian metric of $M$ and $B \in T_B\mathcal{B}^L_{k_0,\dots,k_{L-1}}(U)$. The integral-based metric on $\mathcal{B}^L_{k_0,\dots,k_{L-1}}(U)$ is defined by
		\begin{align*} 
			\langle \langle X, Y \rangle \rangle_{B} := \int_0^L \big\langle X(t), Y(t) \big\rangle_{B(t)} \dd t
		\end{align*}
		for all $X,Y \in T_B\mathcal{B}(U)$.    
	\end{definition} 
	Bézier splines depend smoothly on their control points, hence the structure is smooth.
	The metric gives Bézier splines a natural distance.
	For general manifolds, explicit formulas for geodesics or the Riemannian exponential $\exp^{\mathcal{B}}_B$ and logarithm $\log^{\mathcal{B}}_B$ at $B \in \mathcal{B}(U)$ seem hard to find; nevertheless, one can approximate them efficiently using variational time-discretization (see Sec.~\ref{sec:computation_hierarchical}).
	
	\begin{remark}
		Consider the Bézierfolds discussed in Remark~\ref{rem:constant_bezierfolds} that consist of constant maps to $U$ only. Since their tangent vectors are constant maps into $TU$, they are isometric to $U$ under the map (\ref{eq:diffeo_F}) when endowed with the integral-based metric~(\ref{eq:metric_bezierfold}).
	\end{remark}
	
	\subsubsection{Sasakian Metric for Cubic Splines}\label{sec:sasakian_metric4splines}
	As first done in~\cite{NavaYazdaniAmbellanHaniketal2023}, the (product) Sasaki metric (\ref{def:Sasaki_metric}) can be bulled back to a cubic B\'ezierfold. 
	
	They endow $\mathcal{B}^L_{3,\dots,3}(U)$ with the pullback $\langle \langle \fcdot, \fcdot \rangle \rangle^S$ of the product Sasaki metric~(\ref{def:Sasaki_metric}) under $F$.
	\begin{definition} \label{def:Sasakian_metric}
		Let $X, Y \in T_B\mathcal{B}^L_{3,\dots,3}(U)$ and $(v^X_0, w^X_0), \dots, (v^X_L, w^X_L) \in (TM)^2$ and $(v^Y_0, w^Y_0), \dots, (v^Y_L, w^Y_L) \in (TM)^2$ the vectors that induce them.
		The Sasakian metric on $\mathcal{B}^L_{3,\dots,3}(U)$ is defined by
		\begin{equation*} 
			\langle \langle X, Y \rangle \rangle_B^S := \sum_{i=0}^L \langle v^X_i,v^Y_i \rangle_{p_i} + \langle w^X_i, w^Y_i \rangle_{p_i}.
		\end{equation*}
	\end{definition}
	
	An advantage when using this metric is that $\mathcal{B}^L_{3,\dots,3}(U)$ and $(TM)^{L+1}$ become isometric (with isometry $\widetilde{F}$).
	Hence, all computations can be done in $(TM)^{L+1}$. In particular, it is never necessary to compute vector fields along Bézier splines explicitly; only the vectors that induce it are needed. Still, iterative procedures are often needed for computations in (products of) $TM$; algorithms for the exponential and logarithm maps have been presented in~\cite{MuFl2012, NavaYazdani2022sasaki} and are summarized in Sec.~\ref{sec:computation_hierarchical}.
	
	\section{Regression with Bézier Splines}
	\label{sec:spline_regression}
	Regression techniques are essential in contemporary statistics, as they analyze the relationship between a dependent variable and explanatory variables. Normal multivariate regression methods are inapplicable when the observed variable takes values in general manifolds. Therefore, a lot of effort has been put into generalizing regression to nonlinear spaces. Non-parametric~\cite{DavisFletcherJoshi2007, Yuan_ea2012, MallastoFeragen2018}, semi-parametric~\cite{Shi_ea2009, Zhu_ea2009}, and parametric~\cite{Hinkle2014, Kim_ea2014, Cornea_ea2017} regression models have been studied using Riemannian tools; some of these are applicable to data from general manifolds~\cite{DavisFletcherJoshi2007, Shi_ea2009, Hinkle2014, Kim_ea2014, MallastoFeragen2018} while others are formulated for specific ones~\cite{Zhu_ea2009, Cornea_ea2017, Yuan_ea2012}. 
	
	Parametric models are usually faster to compute and easier to interpret; therefore, they are often preferred when one can be confident that the relationship between the measured and explanatory variables is of a known type. If the latter cannot be guessed, then non-parametric models should be used.
	The simplest parametric regression model in Euclidean space is linear regression, and geodesic regression~\cite{niethammer2011geodesic, Fletcher2013} was introduced as its generalization to Riemannian manifolds. 
	
	Although geodesics are a helpful tool for describing many processes, they are not always accurate for certain types of relationships, such as periodic or saturated processes.
	In such cases, higher-order (parametric) models are adequate. Therefore, Riemannian polynomials~\cite{Hinkle2014} have been considered as underlying curves. They are far more flexible than geodesics and, thus, can describe a much bigger class of relationships. Nevertheless, they often cannot be used due to their high computational complexity. They are defined (through variational principles) as solutions to differential equations involving curvature terms. Since closed-form solutions of the equations are usually not available, evaluation and optimization are complicated and computationally expensive.
	
	The authors of~\cite{Hanik_ea2020} suggest using manifold-valued Bézier splines as an alternative to Riemannian polynomials. Bézier splines have explicit formulas in many manifolds, which makes computations faster. They are also easily concatenated, which is a notable advantage: Riemannian polynomials allow for piecewise composition when degrees are odd, but for even degrees, it is unclear how to define splines~\cite[p.\ 34]{Hinkle2014}. 
	
	In this section, we will provide a summary of how to use Bézier splines for regression in manifolds. It is possible to use either general or closed splines; we formulate everything using the former. For closed splines one must simply replace $\mathcal{B}(U)$ with $\mathring{\mathcal{B}}(U)$ and $K$ with $\mathring{K}$. After defining the model, we discuss the estimation of the involved notions. As a novel contribution that is not included in~\cite{Hanik_ea2020}, we prove the equivalence of the least-squares and maximum-likelihood estimation in symmetric spaces under the assumption of a Gaussian-like distribution for the errors. Finally, we discuss an application of the regression to normalize manifold-valued data w.r.t.\ a confounding variable. This method was introduced in~\cite{Hanik_ea2023}.

	\subsection{The Model}
	For $k$-th order polynomial regression in Euclidean space $\mathbb{R}^d$ one assumes that the relationship between an $\mathbb{R}^d$-valued random variable $Y$ and a scalar variable $t$ is given by the model
	\begin{equation} \label{eq:regression_model_multivariate}
		Y(t) = P(t) + \epsilon(t),
	\end{equation}
	with $P$ being a (deterministic) polynomial of degree at most $k$ and $\epsilon$ a vector-valued, random variable representing the error. In the following, we generalize this model to Riemannian manifolds using Bézier splines; furthermore, we extend it conceptually by including closed curves.
	
\begin{figure}[ht]
		\centering
		\includegraphics[width=.6\linewidth]{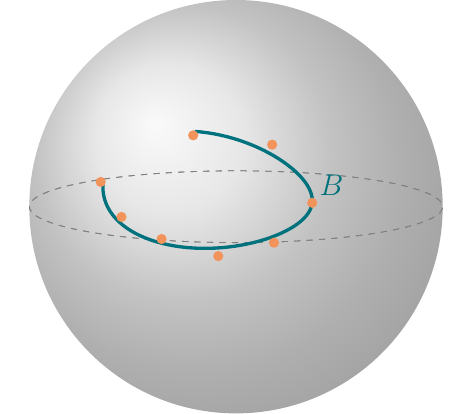}
		\caption{Sketch of regression with non-geodesic Bézier splines in $\mathcal{S}^2$. The orange dots are data points, while $B$ indicates the underlying spline $B$ in~(\ref{eq:regression_model_manifold}).}
		\label{fig:regressed_spline}
	\end{figure}

	Let $M$ be a Riemannian manifold and $U \subseteq M$ be a normal convex neighborhood. Further, let $n$ data points $(t_1,q_1),\dots,(t_n,q_n) \in [0,1] \times U$ with corresponding scalar parameters (e.g., points in time) be given. We suppose that there is a Bézier spline $B \in \mathcal{B}(U)$ such that the data points $q_1,\dots,q_n$ are independent realizations of an $M$-valued random variable $\mathcal{Q}$ that is connected to a deterministic variable $t \in \mathbb{R}$ according to the model
	\begin{equation} \label{eq:regression_model_manifold}
		\mathcal{Q}(t) = \exp_{B(t;p_0,\dots,p_K)} \big(\epsilon(t) \big).
	\end{equation}
	Here, for each $t \in [0,1]$ the random error variable $\epsilon$ takes values in $T_{B(t)}M$ whose realizations are (generally non-continuous) random vector fields along $B$. The assumption that the explanatory variable $t$ is in the interval $[0,1]$ is purely for reasons of consistency: For a finite set of samples, the realizations of $t$ must be in an interval of finite length, which can always be linearly re-scaled to $[0,1]$. 
	The unknown parameters of the model are the control points $p_0,\dots,p_K$. A visualization of the model for samples on the sphere is shown in Fig.~\ref{fig:regressed_spline}.
	
	In Euclidean space, (\ref{eq:regression_model_manifold}) reduces to (\ref{eq:regression_model_multivariate}) when restricting to Bézier curves; i.e., (\ref{eq:regression_model_manifold}) is a generalization of polynomial regression. It also generalizes geodesic regression~\cite[Eqn.\ 3]{Fletcher2013}, to which it reduces when $B$ consists of a single segment with two control points.
	
	In the following, we always assume that the type of $B$ is known. In practice, this must be determined first. As for polynomial regression, knowledge-driven (first) guesses of the degrees should be tested, e.g., by holding back a testing set from the data to validate the model.
	
	\subsection{Least Squares Estimation} \label{sec:least_squares_estimation}
	We now discuss how the Bézier spline $B$ in (\ref{eq:regression_model_manifold}) should be estimated. It is often difficult to use maximum likelihood estimation or empirical Bayes due to the rarity of tractable noise distributions on manifolds. Therefore, geometric least squares estimation, which will be discussed first, is commonly employed. Nevertheless, we will show towards the end of this section that the former coincides with maximum likelihood estimation in some situations.
	
	To motivate our method, we again look at the multivariate case first.
	A standard way to estimate the polynomial $P$ in (\ref{eq:regression_model_multivariate}) is \textit{ordinary least squares}. Let $N$ measurements $(t_1,y_1),\dots,(t_n,y_n) \in [0,1] \times \mathbb{R}^d$ be given.
	Denoting the set of polynomial curves over $[0,1]$ of degree at most $k$ by $\mathcal{P}_k$, we can approximate $P$ by the \textit{least squares estimator}
	$$\widehat{P} := \argmin_{P \in \mathcal{P}_k} \frac{1}{2}\sum_{i=1}^n \|P(t_i) - y_i\|^2.$$
	
	This idea directly translates to Riemannian manifolds. Let $\mathcal{C} \subset U^{K+1}$ be the largest (i.e., containing all others with the same property) \emph{compact}, simply-connected subset such that each point in $\mathcal{C}$ defines an element of $\mathcal{B}(U)$.
	Then, we can define the least squares estimator of (\ref{eq:regression_model_manifold}) as follows.
	\begin{definition} \label{def:least_squares_estimation_riemannian}
		Let $n$ realizations
		$(t_1,q_1),\dots,(t_n,q_n) \in [0,1] \times U$
		of (\ref{eq:regression_model_manifold}) be given. Then, the \textnormal{sum-of-squared error} of a Bézier spline $B \in \mathcal{B}(U)$ is defined by
		\begin{equation} \label{eq: error}
			\mathcal{E}(p_0,\dots,p_K) := \frac{1}{2} \sum^n_{i=1} \dist \Big(B(t_i;p_0,\dots,p_K),q_i \Big)^2.
		\end{equation}
		Furthermore, the \textnormal{least squares estimator} $\widehat{B} \in \mathcal{B}(U)$ of the spline $B$ in (\ref{eq:regression_model_manifold}) is determined by its control points $$(\widehat{p}_0,\dots,\widehat{p}_K) \in \mathcal{C}$$ 
		according to
		\begin{equation} \label{eq:minimizer}
			(\widehat{p}_0,\dots,\widehat{p}_K) := \argmin_{(p_0,\dots,p_K) \in \mathcal{C}}\ \mathcal{E}(p_0,\dots,p_K).
		\end{equation}
	\end{definition}
	Note that, for two control points, this definition also boils down to Fletcher's definition~\cite[Eqn.\ 5]{Fletcher2013} so that in the geodesic case, the least squares estimates from geodesic regression and our method coincide.
	Generally, none of the control points coincides with a data point, as with spline \textit{interpolation}~\cite{BergmanGousenbourger2018, GousenbourgerMassartAbsil2019}.
	
	Before going on, a comment on our use of $\mathcal{C}$ seems in order. Like Fletcher in~\cite[Thms.\ 2 and 3]{Fletcher2013}, we introduce it to ensure the existence of a least squares estimator. The existence question in an open set is very involved and, as far as we can see, depends strongly on the ``layout'' of the data, the geometry of the space, and the degree of the underlying curve. E.g., if we want to use a polynomial of high degree for sparse data with similar values of $t$ but highly varying realizations of $\mathcal{Q}$, this will often lead to strongly oscillating least squares estimators (a well-known effect for polynomials of a higher degree in Euclidean space), the ``extrema'' eventually ``trying'' to leave $U$. We leave this question for future work and use a compact set instead.
	
	In general, the minimizer~(\ref{eq:minimizer}) is not known analytically, which makes iterative schemes necessary. Therefore, we apply Riemannian gradient descent. (For a reference on optimization on manifolds see \cite{AbsilMahonySepulchre2007}.) 
	The gradient of $\mathcal{E}$ can be computed for each control point individually when the product manifold $M^{K+1}$ is endowed with the so-called product structure (see, e.g.,~\cite{doCarmo1992})---a very natural choice since we have no reason to assume that there are dependencies between the control points.
	
	It is shown in~\cite[Sec.\ 4.2]{BergmanGousenbourger2018} that $\grad_{p_j}\mathcal{E}$ can be computed explicitly in symmetric spaces.

	As an initial guess for the gradient descent, choosing $(p_0,\dots,p_K)$ along the geodesic polygon whose corners interpolate the data points closest to knot points (in terms of $t$) has worked well. Note that gradient descent methods converge to \textit{local} minima. Thus, we implicitly assume that there is no ``badly behaved'' local minimum ``far away'' from the least squares estimator.  
	
	\subsection{Maximum Likelihood Estimation}
	In Euclidean space, the least-squares estimation coincides with maximum likelihood estimation under the assumption that the distribution of the errors is Gaussian. In~\cite[Sec.\ 4.4]{Fletcher2013}, Fletcher showed that this is also true for geodesic regression in symmetric spaces under a generalized isotropic normal distribution. 
	In the following, we extend his result to regression with Bézier splines. To this end, we deal with integrals of functions on manifolds; see~\cite[Ch.\ 2]{Jost2017} for a reference.
	
	The maximum likelihood method is another probabilistically motivated method to obtain an estimator for the parameters of a statistical model from given data. The idea is to view the parameters (the control points of $B$ in our case) as variables of the underlying probability density function (pdf). The \textit{maximium likelihood estimator} is the set of parameters that maximizes the so-called \textit{likelihood function} under the given data. In practice, because it has the same extrema but is often easier to analyze, one usually applies the logarithm to the likelihood function and works with the \textit{log-likelihood function} instead. 
	
	An assumption on the pdf underlying the model is needed to apply the maximum likelihood method. Let $M$ be a connected symmetric space, $\overline{q} \in M$, and define
	\begin{equation*}
		C(\sigma) := \int_M \exp \left( -\frac{ \dist(p, \overline{q})^2}{2 \sigma^2} \right) \dd p.
	\end{equation*}
	It is important to note that $C$ does not depend on $\overline{q}$: The distance function and thus $C$ are invariant under isometries, and for any two points in a symmetric space, there always is an isometry of the space mapping one into the other~\cite{Helgason2001}.
	Given $\overline{q}$ and $\sigma > 0$, Fletcher's generalized normal distribution is now defined via the pdf
	\begin{align} \label{eq:isotropic_gaussian}
		\rho&: M \to \mathbb{R}_{\ge 0} \nonumber \\
		p \mapsto \rho(p; \overline{q}, \sigma) &:= \frac{1}{C(\sigma)} \exp \left( -\frac{ \dist(p, \overline{q})^2}{2 \sigma^2} \right).
	\end{align}
	The distance function in the exponential leads to an isotropic assignment of probability mass.
	Two advantageous properties of $\rho$ are (apart from being an isotropic normal distribution in Euclidean space) that it is smooth and nonzero everywhere; both need not hold for other generalizations of the normal distribution. (They are not given, e.g., for the maximum entropy distribution introduced in~\cite{Pennec2006}, which, on the other hand, has the advantage that anisotropy can be modeled.)
	
	We say that the model (\ref{eq:regression_model_manifold}) has \emph{Gaussian errors}, if the random variable $\mathcal{Q}$ is conditionally distributed according to (\ref{eq:isotropic_gaussian}), i.e., with pdf
	\begin{equation}
		\nu(q| t=t_0) = \rho \big(q; B(t_0;p_0,\dots,p_K), \sigma \big)
	\end{equation}
	for some $\sigma > 0$.
	Under this assumption, the \textit{log-likelihood function} of our model is given by
	\begin{align*}
		l\left(p_0,\dots,p_K;\sigma \right) &:= N \log \big(C(\sigma) \big) - \frac{1}{2\sigma^2} \sum_{i=1}^n \dist \big(B(t_i;p_0,\dots,p_K), q_i \big)^2
	\end{align*}
	for given data $(t_1,q_1),\dots,(t_n,q_n) \in [0,1] \times U$;
	the \textit{maximum likelihood estimator} of the model then is the spline whose control points maximize $l$ (with control points in $\mathcal{C} \subseteq U^{K+1}$).
	
	The following theorem shows that, in symmetric spaces, maximum likelihood and least squares estimators coincide under the above assumption.
	
	\begin{theorem} \label{thm:max_like}
		Let $M$ be a symmetric space, $U \subseteq M$ a normal convex neighborhood, and
		$(t_1,q_1), \dots, (t_n,q_n) \in [0,1] \times U$ realizations of the model (\ref{eq:regression_model_manifold}). Then, least squares optimization and maximum likelihood approximation are equivalent for (\ref{eq:regression_model_manifold}) with Gaussian errors.
	\end{theorem}
	\begin{proof}
		\citet[Sec.\ 4.2]{BergmanGousenbourger2018} showed that the gradients 
		$$\grad_p \dist \big( B(t_i;p_0,\dots,p_{j-1},\fcdot,p_{j+1},\dots,p_K), q_i \big)^2, \quad j=0,\dots,K,$$ 
		exist for any $p \in U$ by deriving an explicit procedure for their computation. Therefore, for each $j=1,\dots, K$ and $p \in U$ we find
		\begin{align*}
			\grad_p l &= -\frac{1}{2\sigma^2} \sum_{i=1}^n \grad_p \dist \big( B(t_i;p_0,\dots,p_{j-1},\fcdot,p_{j+1},\dots,p_K), q_i \big)^2 \\ 
			&= -\frac{1}{\sigma^2} \grad_p\mathcal{E}. 
		\end{align*}
		Hence, in the interior of $\mathcal{C}$, local minimizers of the sum-of-squared error $\mathcal{E}$ are local maximizers of the log-likelihood function $l$, and vice versa. 
		
		We now show that this also holds for global extrema in $\mathcal{C}$. Let $(\widehat{p}_0,\dots,\widehat{p}_K)$ be the global minimizer of $\mathcal{E}$, $(p_0,\dots,p_K) \in \mathcal{C}$ arbitrary, and $\alpha$ a smooth curve in $\mathcal{C}$ that connects them, while ending in the minimizer. Using the definition~(\ref{eq:gradient}) of the gradient and the relationship above, we see that the differentials of $l$ and $\mathcal{E}$ are also related by $\dd l = -1/\sigma^2\ \dd \mathcal{E}$. Hence, using integrals of differential 1-forms~\cite[Sec.\ 2]{Jost2017}, we can apply Stoke's Theorem~\cite[Thm.\ 2.1.6]{Jost2017}:\footnote{Stoke's theorem acts as the manifold-version of the fundamental theorem of calculus for line integrals here.} 
		\begin{align*}
			l(\widehat{p}_0,\dots,\widehat{p}_K) - l(p_0,\dots,p_K) &= \int_\alpha \dd l \\ 
			&= -\frac{1}{\sigma^2} \int_\alpha \dd \mathcal{E} \\
			&= \frac{1}{\sigma^2} \big(\mathcal{E}(p_0,\dots,p_K) - \mathcal{E}(\widehat{p}_0,\dots,\widehat{p}_K)\big) \\
			&\ge 0.
		\end{align*}
		Since $(p_0,\dots,p_K)$ is arbitrary, this shows that $(\widehat{p}_0,\dots,\widehat{p}_K)$ is a global maximum of $l$. Analogously, one can show that global maxima of $l$ are global minima of $\mathcal{E}$.
	\end{proof}

	\subsection{Normalization via Regression} \label{sec:normalization}
	
	In statistics, one almost always needs to control confounding and other extraneous variables that are not of interest to a study but might influence the outcome. 
	If the pool of samples is large enough, the controlling can be done through sample selection, enforcing, e.g., that the variable is constant or follows a particular distribution. Another possibility, which is also applicable in the case of a small data pool, is \textit{normalization}, i.e., shifting/scaling the data's statistics such that the influence of the confounding variables is minimized.
	Such procedures can also be necessary when analyzing manifold-valued data. In this section, we summarize the Bézier spline-based normalization method from~\cite{Hanik_ea2023} that helps to reduce the interference of continuous influence parameters. 
	While motivated by the a specific application, it can be applied when several sets of manifold-valued data shall be compared that are influenced by a continuous confounding variable. 
	Conceptually similar procedures are used in Euclidean space, e.g., when analyzing molecules~\cite{Shen_ea2016, Wahid_ea2016, HafenmeisterSatija2019}. 
	
	\subsubsection{The Model}
	Let $M$ be a Riemannian manifold, $U \subseteq M$ a normal convex neighborhood, and $I_s \subset \mathbb{R}$, $s=1,\dots,S,$ closed intervals with non-empty intersection, i.e.,\ $\cap_{s=1}^S I_s \ne \emptyset$. 
	(The latter assumption guarantees a \textit{shared} parameter value that can serve as a reference point for the normalization.)
	Let further $s$ groups of $U$-valued data points be given, each element coming with a parameter: 
	\begin{equation*}
		\left( t^{(s)}_j, q^{(s)}_j \right) \in I_s \times U, \quad s=1,\dots,S, \quad j=1,\dots,n_s.
	\end{equation*}
	We assume that the data from $s$-th group behaves according to the model (\ref{eq:regression_model_manifold}), i.e.,
	\begin{equation} \label{eq:model_normalization}
		\mathcal{Q}^{(s)}(t) := \exp_{B^{(s)}(t;p_0,\dots,p_K)} \left(\mathcal{R}^{(s)}(t) \right) \quad s=1,\dots,S.
	\end{equation}
	In contrast to the previous sections, here it is supposed that we are interested in the $TM$-valued variables $\mathcal{R}^{(s)}$ (thus the change of notation). We further assume that the ``drifts`` in the data given by the Bézier splines $B^{(s)}$, $i=1,\dots, S$, are \textit{caused by the single confounding but deterministic variable $t$}. 
	
	\subsubsection{Normalization}
	In this section, we discuss how the data can be normalized at some $t_0 \in \cap_{i=1}^S I_i$ so that further analysis is not influenced by the variability caused by $t$. The procedure is as follows.
	First calculating the least squares estimators (\ref{eq:minimizer}) for each group yields best fitting splines $\widehat{B}^{(s)}: I_i \to M$, $i = 1,\dots,S$. 
	It is possible that the regressed curves already provide valuable information about the data, as we will see in Sec.~\ref{sec:archaeology}.
	Approximations $\widehat{R}^{(s)}(t_i) \in TM$ of the realizations of $\mathcal{R}^{(s)}(t_i) \in TM$ of (\ref{eq:model_normalization}) are then given by the logarithms
	\begin{equation*}
		\widehat{R}^{(s)}(t_i) := \log_{\widehat{B}^{(s)}(t_j^{(s)})} \left( q^{(s)}_j \right) \in T_{\widehat{B}^{(s)}(t_j^{(s)})}M, \quad i=1,\dots,S, \quad j=1,\dots,n_i.
	\end{equation*}
	To minimize the influence caused by the variable $t$, we want to normalize the data at some chosen point $t_0$.
	Since, for each group, $\widehat{B}^{(s)}$ approximates the trend that $t$ enforces, we propose to 
	parallel translate the vectors $\widehat{R}^{(s)}_j$ along $\widehat{B}^{(s)}$ to $\widehat{B}^{(s)}(t_0)$.
	This results in vectors
	\begin{equation} \label{eq:linear_normalized_data}
		w^{(s)}_j \in T_{\widehat{B}^{(s)}(t_0)}M, \quad i=1,\dots,S, \quad j=1,\dots,n_i.
	\end{equation}
	They represent the differences of the data points to the Bézier splines \textit{normalized at} $t_0$. Mapping them back to the manifold gives the \textit{normalized data points}:
	\begin{equation} \label{eq:normalized_data}
		\widetilde{q}^{(s)}_j := \exp_{\widehat{B}^{(s)}(t_0)} \left( w^{(s)}_j \right), \quad i=1,\dots,S, \quad j=1,\dots,n_i.
	\end{equation}
	\begin{figure*}[ht!]
		\centering
		\includegraphics[trim=5cm 4cm 5cm 6cm, clip,width=.85\textwidth]{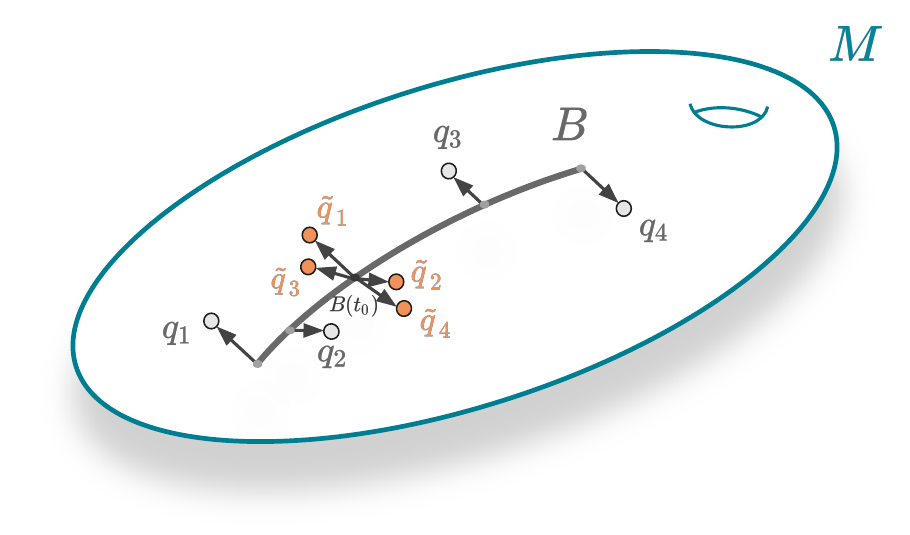}
		\vspace{-2em}
		\caption{Normalization w.r.t.\ some parameter $t$ for a single data group $(q_j,t_j)$, $j=1,2,3,4$, in a Riemannian manifold $M$. The curve $B$ is the result of spline regression w.r.t.\ $t$. The points $B(t_j)$ are depicted in light grey, while the tangent vectors $v_j$ are the black arrows attached to them. Finally, the parallel translation $w_j$ of each $v_j$ is a tangent vector at $B(t_0)$ (also black); it yields the normalized data $\widetilde{q}_j$ shown in orange.}\label{fig:normalization}
		\hfill
		\centering
		\includegraphics[trim=5cm 4cm 5cm 6cm, clip,width=.85\textwidth]{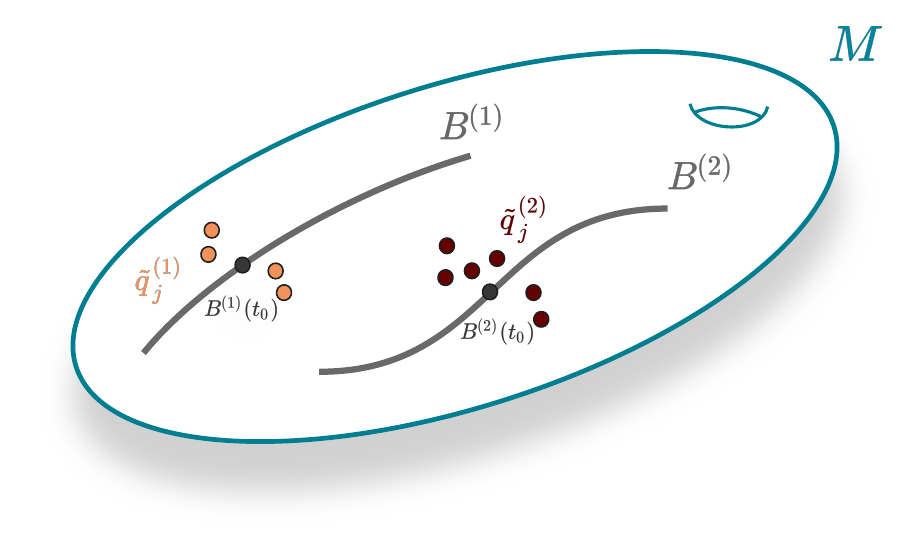}
		\vspace{-2em}
		\caption{Normalization results for two groups. Both splines $B_1$ and $B_2$ are defined for the value $t = t_0$, so the two data sets can be normalized at $B^{(1)}(t_0)$ and $B^{(2)}(t_0)$, respectively. The groups (orange and red) can now be compared without bias caused by the influence of $t$.}\label{fig:comparison}
	\end{figure*}
	

	Now, we can perform \emph{inter-group comparison} with the normalized data (\ref{eq:normalized_data}). E.g., we can compute the Fréchet means~\cite[Ch.\ 2]{PennecSommerFletcher2020} and perform group tests for equality as in~\cite[Sec.\ 3.3]{MuFl2012}. If the curvature near the points $\widehat{B}^{(s)}(t_0)$ is small (i.e., all sectional curvatures are close to $0$), then we can also use the linearized data (\ref{eq:linear_normalized_data}) and apply methods from multivariate statistics; the higher the curvature, though, the stronger will be the introduced distortion.
	The normalization process and its result for data from two groups are visualized in Figs.~\ref{fig:normalization} and \ref{fig:comparison}.
	It is important to note that the choice of $t_0$ will influence the results, just like in multivariate statistics. (Imagine, e.g., that we measure the height of a group of people. Only considering children at age ten will give results different from 20-year-old adults.)
	
	Computation-wise, for geodesics as underlying trends, explicit formulas for parallel transport are known for many manifolds that appear in applications. Otherwise, several numerical schemes can provide approximations of parallel transport along curves~\cite{KheyfetsMillerNewton2000, Louis_ea2018, GuiguiPennec2021}.

	\section{Hierarchical Models} \label{sec:hierarchical_model}
	
\begin{figure*}[bt]
		\centering
		\includegraphics[trim=5cm 4cm 5cm 6cm, clip,width=.85\textwidth]{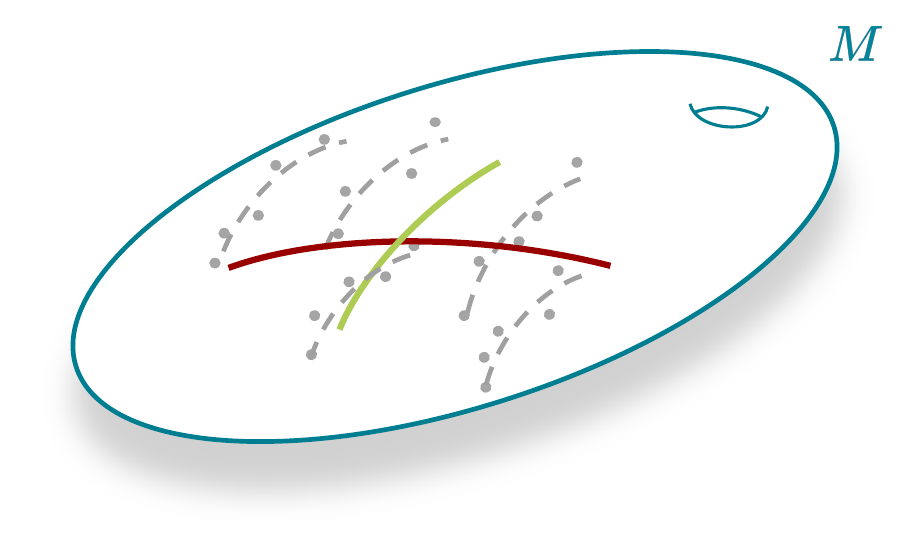}
		\vspace{-2em}
		\caption{Conceptual comparison of a hierarchical model and regression analysis on a longitudinal data set. The measurements are shown as grey dots (without time). Each subject's (correlated) measurements follow its trend (broken lines). Spline regression on the full data set yields the red curve; its direction deviates strongly from all individual trends. The mean of the individual trends (green curve) is a better estimator of the average longitudinal trend in the cohort.}
		\label{fig:hierarchical_model}
	\end{figure*}

	Time-dependent data analysis is becoming increasingly relevant for a wide range of applications, including the investigation of disease onset and progression, physical performance assessment from biomechanical gait data, and face expression analysis in video sequences.
	All of these are examples of longitudinal data, in which individual instances of a common underlying process are observed repeatedly over a period of time.
	While popular statistical tools such as mean-variance analysis and regression allow for the investigation of phenomena among individuals or inside a single one, longitudinal data contains correlations within measurements of single individuals, which violates the independence assumptions of such cross-sectional methods.
	Figure~\ref{fig:hierarchical_model} sketches the difference between cross-sectional and longitudinal models.
	Another concern is missing data, which can occur as a result of acquisition problems or when subjects drop out of a clinical trial.
	Proper statistical inference for longitudinal data, thus, needs to account for both within-individual correlations of observations as well as the potentially sparse or non-uniform sampling.
	
	From this perspective, hierarchical models---specifically mixed-effects models---offer a suitable and versatile framework for the analysis of longitudinal data~\cite{locascio2011overview, GFS2016}.
	These approaches address the multitude of inherent interrelations by formulating a model that assumes that each subject possesses an individual functional relationship between the dependent variable and the time-related predictor(s).
	Therefore, a parametric spatiotemporal model is estimated to fit the data for each individual.
	Any discrepancy in the fit is attributed to a random error variation in the observed dependent variable.
	The subject-specific coefficients that characterize these models (such as the intercept and slope in the case of straight lines) are assumed to exhibit random variation within the population.
	In addition to these random effects, there exists a ``fixed'' effect that is frequently the main focus of interest.
	This effect is defined as a ``fixed'' coefficient that represents the average spatiotemporal model for the entire group.
	
	The existing body of literature extensively covers the mixed-effects framework for scalar and vector-valued measurements. However, the extension of this framework to other types of data with structured domains, such as shapes or graphs, has only emerged in recent years.
	Based on geodesics, i.e., generalizations of straight lines, multiple works~\cite{MuFl2012, singh2013hierarchical, Hong_ea2015, muralidharan2017diss, NavaYazdaniHegevonTycowicz2022} derived geodesic hierarchical models that encode both subject-specific as well as group trends in terms of geodesics.
	Nonmonotonous shape changes, e.g., present in a time series of cardiac motion or anatomical changes in the human brain over decades, do not, in general, adhere to constraints of geodesicity.
	This assumption can, hence, be a limiting factor diminishing the flexibility and fidelity of such approaches, calling for the development of higher-order models.
	
	\subsection{The Model}
	
	In line with the regression approach presented in the previous section, we can obtain higher-order, hierarchical models building upon intrinsic B\'ezier splines. 
	We now introduce such a nonlinear hierarchical model following the derivation in~\cite{HanikHegevonTycowicz2022}.
	Since a key application of mixed-effects models is in life sciences, we speak of measurements taken from ``subjects'' in the following.
	Nevertheless, the model is not restricted to data from this field.

	Let $M$ be a Riemannian manifold and $U \subseteq M$ a normal convex neighborhood. Further, let $\mathcal{B}(U)$ be a Bézierfold over $U$.
	Consider that $S$ subjects are given and that, for each, there is (a possibly different number of) $n_s$, $s=1,\dots, S$, measurements of the pair of an independent, deterministic scalar variable $t$ and a manifold-valued dependent variable $\mathcal{Q}^{(s)}$. The data thus reads 
	\begin{equation} \label{eq:longitudinal_data}
		\left( t_i^{(s)}, q_i^{(s)} \right) \in \mathbb{R} \times U, \quad  i=1,\dots,n_s, \quad s=1,\dots,S.
	\end{equation}
	%
	While cross-sectional measurements (i.e., points $q_i^{(s)}$ with mutually different $s$) can be assumed as being independent of each other, intra-subject measurements (i.e., samples $q_i^{(s)}$ with the same $s$) are correlated and, thus, have to be treated separately.
	To this end, a hierarchical two-stage model consisting of the following two consecutive levels can be employed.
	
	\noindent\textbf{Individual level}: On the lower, individual level, we assume that (\ref{eq:regression_model_manifold}) underlies each subject's data independently from the others; the type and the degrees of all Bézier splines are supposed to be the same; i.e., it is assumed that the data (\ref{eq:longitudinal_data}) is drawn from random variables
	\begin{equation*} \label{eq:model_hierarchical_1}
		\mathcal{Q}^{(s)}(t) := \exp_{B^{(s)}(t;p_0^{(s)},\dots,p_K^{(s)})} \left(\epsilon^{(s)}(t) \right), \quad s=1,\dots,S.
	\end{equation*}
	\noindent \textbf{Group level:} On the upper group level, we view each trajectory as a perturbation of a common mean trajectory $\overline{B}:= \overline{B}(\fcdot;\overline{p}_0,\dots,\overline{p}_K) \in \mathcal{B}(U)$ in the Bézierfold $\mathcal{B}(U)$ according to
	\begin{equation*} \label{eq:model_hierarchical_2}
		B^{(s)} = \exp^{\mathcal{B}}_{\overline{B}} \left(X^{(s)} \right), \quad s=1,\dots,S,
	\end{equation*}
	with random tangent vector fields $X^{(s)} \in T_{\overline{B}}\mathcal{B}(U)$.
	This model first deals with the correlated samples encoding them as a single, subject-specific spline model, so, on the group level, independence can be assumed again. 
	
	\subsection{Parameter Estimation}\label{sec:computation_hierarchical}
	As the individual level does not depend on the group level parameter, the control points of the individual developments and those of the mean trend can be estimated efficiently by a two-step least squares procedure.
	In particular, in the first step, the parameters of the individual trends can be determined in parallel using the estimator (\ref{eq:minimizer}) as discussed in Sec.~\ref{sec:least_squares_estimation}.
	The second step then consists of finding the Fréchet mean $\overline{B}$ in $\mathcal{B}(U)$ of the estimated individual trends.
	The following discussion is, therefore, only concerned with approximating the control points of the mean spline.
	As the Fr\'echet mean~(\ref{def:frechet_mean}) is defined with respect to geodesic distances, the estimated model depends on a choice of Riemannian metric for $\mathcal{B}(U)$.
	
	\subsubsection{Integral-based Metric}
	For the Riemannian metric discussed in Sec.~\ref{sec:integral-based_metric}, we can obtain an estimation scheme in terms of the discrete geodesic calculus by \citet{RumpfWirth2014} that employs a variational time-discretization.
	
	Let $B_1, B_2 \in \mathcal{B}^L_{k_0,\dots,k_{L-1}}(U)$. A path between $B_1$ and $B_2$ through $\mathcal{B}^L_{k_0,\dots,k_{L-1}}(U)$ (and defined in $[0,1]$) may be represented as a parametrized surface in $U$, because it induces a map $H:[0,1] \times [0,L] \to U, (r, t) \mapsto H(r,t)$, with $H(0, \fcdot) = B_1$ and $H(1,\fcdot) = B_2$.
	A geodesic between $B_1$ and $B_2$ is a minimizer of the \emph{path energy}\footnote{Indeed, a curve minimizes the path energy \textit{if and only if} it is a geodesic; see, e.g.,~\cite[p.\ 196]{doCarmo1992}.}
	\begin{equation*}
		E(H) := \int_0^1 \int_0^L \left\langle \frac{\dd H}{\dd r}(r,t), \frac{\dd H}{\dd r}(r,t) \right\rangle \textnormal{d}t\, \textnormal{d}r.
	\end{equation*}
	
	Discretizing in $\mathcal{B}^L_{k_0,\dots,k_{L-1}}(U)$ and identifying splines with their control points, we obtain a \emph{discrete $\ell$-geodesic} $(p_0^j,\dots,p_K^j)_{j=0,\dots,\ell} \in (U^{K+1})^{\ell + 1}$ between $B_1$ and $B_2$ as the minimizer of the \emph{discrete path energy}
	\begin{align*} \label{eq:E_n}
		\begin{split}
			E_\ell \bigg( \left(p_0^j,\dots,p_K^j \right)_{j=0,\dots,\ell} \bigg) := \ell \sum_{j=1}^{\ell-1} \int_0^L \dist \bigg( B \Big(t; p_0^j,\dots, &p_K^j \Big), \\ 
			&B \Big(t; p_0^{j+1},\dots,p_K^{j+1} \Big) \bigg)^2 \textnormal{d}t.
		\end{split}
	\end{align*}
	When distances on $M$ can be computed, the integral can be evaluated using a suitable quadrature rule. 
	
	To approximate minimizers of $E_\ell$, \citet{HanikHegevonTycowicz2022} extend the iterative procedure from~\cite{NavaYazdaniHegevonTycowicz2022} to this setting. 
	The algorithm is motivated by the general characteristic of shortest paths that are induced by Srivastava's metric and the corresponding connection in the space of curves: If $\alpha_1, \alpha_2 : [0, L] \to U$ are two smooth curves, then the map $H(\fcdot, t): [0,1] \to U$ that is induced by the shortest curve between $\alpha_1$ and $\alpha_2$ is a geodesic in $M$ for all $t \in [0, L]$; see~\cite{SrivastavaKlassen2016}.
	Approximating the integrals in $E_\ell$ with a quadrature and applying an alternating optimization scheme (i.e., block coordinate descent) lets us compute the discrete $\ell$-geodesics between two curves by iteratively performing spline regression.
	First, we initialize the control points of the inner curves equidistantly along the geodesics that connect the corresponding control points of $B_1$ and $B_2$. The inner trajectories are then updated so that they lie ``in the middle'' of their neighbors; to this end, we replace them with the result of a spline regression on the $K+1$ data points given by (equidistant) evaluations of the neighboring curves.

	Next, we discuss the computation of the \emph{discrete $\ell$-mean} of $S$ curves $B_1,\dots, B_S \in \mathcal{B}^L_{k_0,\dots,k_{L-1}}(U)$, with which we approximate the common mean (curve). It is the spline $\overline{B} \in \mathcal{B}^L_{k_0,\dots,k_{L-1}}(U)$ minimizing
	\begin{align*}
		G_\ell(p_0,\dots,p_K) &:= \sum_{s=1}^S\ E_\ell\bigg( \left(p_0^j,\dots,p_K^j \right)_{j=0,\dots,\ell}^{(s)} \bigg), \\
		\textnormal{s.t.\ } (p_0^\ell,\dots,p_K^\ell)^{(s)} &= (p_0,\dots,p_K)^{(s)}, \quad s=1,\dots,S,
	\end{align*}
	where $(p_0^j,\dots,p_K^j)_{j=0,\dots,\ell}^{(s)}$ denotes the control points of the discrete $\ell$-geodesic between $B_s$ and $B(\fcdot;p_0,\dots,p_K)$.
	It can be computed with an alternating optimization scheme. 
	We initialize the control points of $\overline{B}$ with the Fréchet means of the corresponding control points of the data curves. Then, we compute discrete geodesics toward the mean and update the latter by a spline regression. This process is repeated in an alternating fashion.

	\subsubsection{Sasakian Metric}
	To estimate the Fréchet mean~(\ref{def:frechet_mean}) in a cubic Bézierfold $\mathcal{B}^L_{3,\dots,3}(U)$ endowed with the Sasakian metric given in Def.~\ref{def:Sasakian_metric}, we can employ the common Newton-type iteration 
	$$ \overline{B}^{i+1} = \exp^{\mathcal{B}}_{\overline{B}^i} \left( \frac{1}{S} \sum_{s=1}^S \log^{\mathcal{B}}_{\overline{B}^i} (B_s) \right). $$
	Since all computations in $\mathcal{B}^L_{3,\dots,3}(U)$ can be performed in $(TM)^{L+1}$ (cf.~Sec.~\ref{sec:sasakian_metric4splines}), it is sufficient to have procedures to compute the exponential and logarithmic maps on $TM$ under the Sasaki metric.
	As the former computes geodesics on $TM$ for a given pair of initial position $(p_0, u_0) \in TM$ and velocity $(v_0, w_0) \in T_{(p_0,u_0)}TM$, it can be approximated by shooting a geodesic~\cite{MuFl2012} using a forward Euler integration of the geodesic equations (\ref{eq:sasaki_geodesic}) given by
	\begin{align*}
		p_{k+1} &= \exp_{p_k}(\varepsilon v_k), \\
		u_{k+1} &= \phi(u_k + \varepsilon w_k, p_{k+1}), \\
		v_{k+1} &= \phi(v_k - \varepsilon R(u_k,w_k)v_k, p_{k+1}), \\
		w_{k+1} &= \phi(w_k, p_{k+1}),
	\end{align*}
	with step size $\varepsilon = \nicefrac{1}{\ell}$ and map $\phi(v,p)$ realizing the parallel transport of vector $v$ to point $p$.  
	
	Given two points $(p_0, u_0) \in TM$ and $(p_\ell, u_\ell) \in TM$, the logarithmic map is the inverse of the exponential map, thus, returning the initial velocity $(v_0, w_0) \in T_{(p_0,u_0)}TM$ of the geodesic segment that connects both points.
	As for the integral-based metric, we utilize a time-discrete path denoted by $(p_k, u_k)_{k=0,\ldots,\ell}$ that is relaxed iteratively to a discrete geodesic by minimizing the path energy  
	$$ E^S_\ell \left((p_k, u_k)_{k=0,\ldots,\ell}\right) = \frac{\varepsilon}{2} \sum_{k=1}^\ell \left(||v_k||^2 + ||w_k||^2 \right), $$
	where $v_k, w_k \in T_{p_k}M$ encode the change in $p_k, u_k$ via forward finite differences, respectively.
	Discrete geodesics, i.e., minimizers of $E^S_\ell$, are then computed via (intrinsic) steepest descent into the negative gradient direction, which is in turn given by
	\begin{align*}
		\grad_{p_k}E^S_\ell &= -\varepsilon \left(\nabla_{v_k}v_k + R(u_k, w_k)v_k \right), \\ 
		\grad_{u_k}E^S_\ell &= -\varepsilon\, \nabla_{v_k}w_k.
	\end{align*}
	Concrete expressions in terms of intrinsic, finite differences for the involved quantities can be found in~\cite{NavaYazdani2022sasaki}.
	
	\subsection{Second-order Statistics and Group Tests}
	Capturing higher-order statistics can provide a more comprehensive description of the underlying distributions that go beyond the estimation of group average trends.
	In particular, the presented Riemannian metrics allow for geometric statistical tools such as principal geodesic analysis (PGA)~\cite{fletcher2004pga} to be applied.
	PGA provides the estimation of the variance and principal directions of the distribution of trajectories in $\mathcal{B}(U)$ that summarize variations within the time-dependent data.
	The representation of the longitudinal observations within the hierarchical basis of principal components in turn yields a highly compact, yet descriptive encoding that is amenable for downstream tasks such as classification and visualization.
	
	One of the major motivations of longitudinal data analysis is to test whether differences found between two groups are statistically significant or due to random chance.
	To this end, Riemannian generalizations of two-sample hypothesis tests can be applied.
	A prominent example is the Hotelling $T^2$ statistic, which can be thought of as a squared Mahalanobis distance between group-wise means with respect to a pooled sample covariance matrix.
	As there is no straightforward way for pooling covariance matrices defined in different tangent spaces, \citet{MuFl2012} proposed a manifold version, which averages the differences in means gauged in each mean's tangent space weighted by the respective single-group covariance.
	This statistic can further be extended to also account for differences in second-order statistics yielding a generalized Bhattacharyya distance~\cite{Hong_ea2015}.
	Regardless of the used statistic, it is typically difficult to formulate the corresponding parametric distribution for samples on general manifolds.
	Moreover, even if such a formulation would be at hand, it is often preferable not to impose the underlying assumptions on the data distributions.    
	Both concerns can be circumvented by estimating the distribution of the statistic employing a non-parametric permutation test design.
	
	\section{Applications} \label{sec:applications}
	\subsection{Background: Shape Analysis}
	An object's geometric characteristics that are invariant through similarity transformations, i.e., under translations, rotations, and scalings, are collectively referred to as its ``shape'' in mathematics. Different methods are used to encode this information depending on the application. As a result, many ``shape spaces'' that can be utilized to model object shapes exist.
	Kendall shape space~\cite{Kendall_ea2009}, which uses landmark configurations as its foundation, is a prime example; \cite{DrydenMardia2016} is an excellent textbook on shape spaces based on landmarks. 
	References ~\cite{BauerBruverisMichor2014, SrivastavaKlassen2016, Younes2019} provide further in-depth treatments on shape spaces that rely on the group of diffeomorphisms. The discussion of skeletal models can be found in~\cite{SiddiqiPizer2008, Pizer_ea2020}.
	Last, research on physics-based spaces has been done in~\cite{Heeren2014, vonTycowicz_ea2018, AmbellanZachowTycowicz2021}. 
	
	Full invariance under similarity transformations is not implemented in all of these spaces for a variety of reasons, including methodological ones (such as when size information is considered crucial), to keep computational costs in check, or because the representation does not support all three (though it is still useful in applications). Nevertheless, invariance is always attained to some extent.
	Nearly all of the aforementioned shape spaces have (at least locally) a non-Euclidean manifold structure. As a result, they require methods from geometric statistics for data analysis.
	
	In all the applications of shape analysis that we describe in the following, the shape space from~\cite{vonTycowicz_ea2018} was used; we denote it in the following by $\Sigma$.
	
	\subsection{Shape Analysis in Medicine}
	Several of the methods that were discussed in the previous sections were successfully applied to shape analysis problems in medicine; the importance of the underlying geometry for shape data analysis is well-established in this field~\cite{PennecSommerFletcher2020}. In the following, we summarize the findings.
	
	\subsubsection{Reconstruction of Shape Trajectories of the Mitral Valve} \label{sec:appl_mitral_valve}
	\begin{figure}[t]
		\begin{center}
			\includegraphics[width=\textwidth]{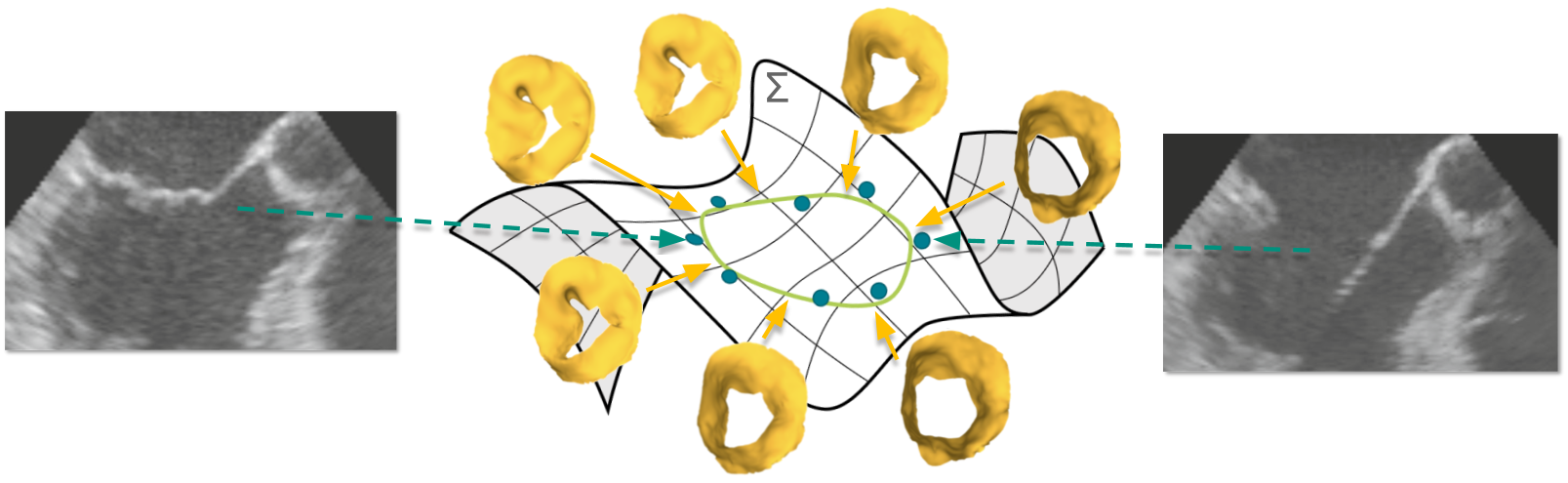}
			\caption{Reconstructed meshes from regression of longitudinal mitral valve data covering a full cardiac cycle. The spline consists of two cubic segments fitted to surfaces extracted from ultrasound images.}
			\label{fig:valve}
		\end{center}
	\end{figure}
	
	Diseases of the mitral valve are the second-most common form of valve disease in adults needing surgery~\cite{Madesis_ea2014}.
	\textit{Mitral valve regurgitation} is an important example thereof. It can have different causes and is characterized by a (backward) blood flow from the left ventricle into the left atrium during systole. This is possible, e.g., when the valve's leaflets do not close fully or prolapse into the left atrium during systole~\cite{EnriquezAkinzVahanian2009}. 
	Mitral valve diseases are often characterized by specific movement patterns. The corresponding shape anomalies can be observed (at least) at some point in the cardiac cycle. Early detection and assessment of mitral valve regurgitation are necessary for the best short-term and long-term results of treatment~\cite{EnriquezAkinzVahanian2009}. 
	Therefore, reconstruction and simulation of 3D mitral valve geometries are active areas of research~\cite{Tautz_ea2020, Walczak_ea2022}. 
	
	In~\cite{Hanik_ea2020}, regression with closed Bézier splines was applied to longitudinal shape data of a mitral valve from a patient with mitral valve regurgitation.
	To this end---sampling the first half of the cycle (closed to fully open) at equidistant time steps---five meshes 
	were extracted from a 3D+t transesophageal echocardiography (TEE) sequence as described in~\cite{Tautz_ea2020}. To reconstruct a full motion cycle, the same five shapes were used in reversed order as data for the second half of the curve. 
	Because of the periodic behavior, a closed spline $\mathcal{\mathring{B}}^2_{3,3}(U)$, $U \subset \Sigma$, with two cubic segments was chosen as the model and an equidistant distribution of the data points along the spline was assumed; i.e., with $q_1,\dots,q_5$ representing the shapes in the space of differential coordinates, the authors employed $$\big(\left(0,q_1\right), \left(\nicefrac{1}{4},q_2\right), \left(\nicefrac{1}{2},q_3\right), \left(\nicefrac{3}{4},q_4\right), \left(1,q_5\right),\left(\nicefrac{5}{4},q_4\right),\left(\nicefrac{3}{2},q_3\right),\left(\nicefrac{7}{4},q_2\right)\big)$$ 
	as the full data set. 
	
	The regressed cardiac trajectory is shown in Fig.~\ref{fig:valve}.
	The spline regression successfully estimates the cyclic motion of the valve and also captures the prolapsing posterior leaflet. The former is remarkable since other regression methods in manifolds do have a hard time modeling cyclic motion (cf.\ Sec.~\ref{sec:spline_regression}), while the latter shows that relevant features are preserved.
	The experiment thus demonstrates the possibility for enhanced reconstruction of mitral valve motion. 
	
	\subsubsection{Remodeling of Knee Bones under Osteoarthritis}
	Osteoarthritis (OA) is a degenerative disease of the joints that affects millions of people worldwide~\cite{ArdenNevitt2006}.
	It develops when the protective cartilage that cushions the ends of the meeting bones degenerates. While the joint pathology is diverse, the most prominent features are the loss of articular cartilage and remodeling of the adjacent bones; the former also leads to the so-called joint space narrowing as less cartilage is there to keep the bones apart. Although any joint can be affected, the population impact is highest for OA of the hip and knee~\cite{ArdenNevitt2006}.
	While OA treatment traditionally consists of pain management and joint replacements for patients with severe symptoms, an improved understanding of the pathogenesis is shifting the focus to disease prevention~\cite{Glynjones_ea2015}.
	\bigskip
	
	\noindent \textit{Regression for OA severity}: In~\cite{Hanik_ea2020}, the authors also use regression with Bézier splines to investigate the development of the shape of the distal femur under OA. Therefore, they regressed the 3D shape against OA severity as determined by the KL grade~\cite{KellgrenLawrence1957}: an ordinal scoring system with grades 0 (healthy) to 4 (severe OA).
	The data set comprised 100 shapes (20 per grade) of randomly selected subjects from the Osteoarthritis Initiative\footnote{\url{https://nda.nih.gov/oai/}} from which they extracted triangle meshes of each bone.
	
	Associating the value $t_i = i/4$ with shapes of grade $i$, they computed the best-fitting curves in $\mathcal{B}^1_1(U), \mathcal{B}^1_2(U), \mathcal{B}^1_3(U)$, $U \subset \Sigma$, i.e., the least squares estimators for regression with geodesic, quadratic, and cubic Bézier curves. The $R^2$-values~\cite[p.\ 56]{PennecSommerFletcher2020} of the regressed curves strongly indicated that the cubic model has far more explanatory power. Particularly, regression with Bézier curves yields results superior to those of standard geodesic regression. A further test also showed that the difference to the geodesic model is not only due to a reparametrization.
	
	\begin{figure}[t]
		\begin{center}
			\includegraphics[width=\textwidth]{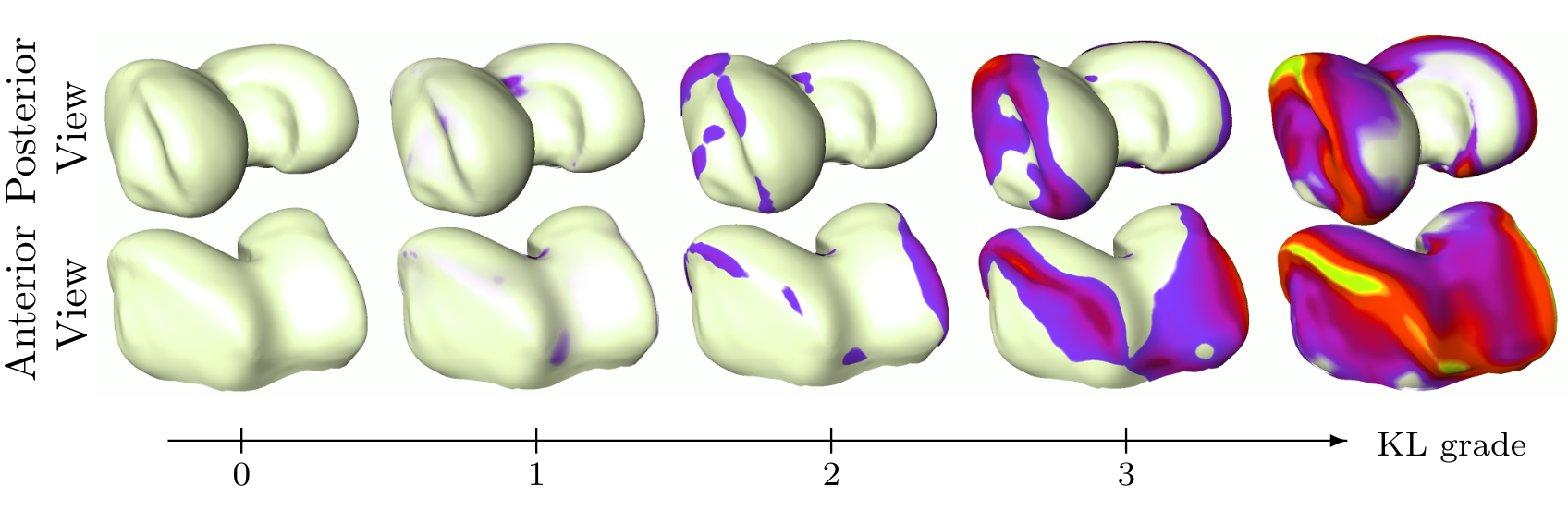}
			\caption{Cubic regression of distal femora.
				Healthy regressed shape ($\textnormal{KL}=0$) together with subsequent grades overlaid wherever the distance is higher than 0.5mm, colored accordingly (0.5 mm~\ShowColormap~3.0 mm).}
			\label{fig:femur_cubic}
		\end{center}
	\end{figure}
	
	The computed cubic Bézier curve is displayed in Fig.~\ref{fig:femur_cubic}. The obtained shape changes consistently describe OA-related malformations of the femur, namely, the widening of the condyles and osteophytic growth. Furthermore, we observe only minute bone remodeling for the first half of the trajectory, while accelerated progression is visible for the second half.
	\bigskip
	
	\noindent \textit{Hierarchical modeling of the femur under OA}: In the last application, the authors investigated the shape development of the femur under OA by using \textit{inter-subject} data (i.e., only one measurement from each subject). 
	In~\cite{HanikHegevonTycowicz2022}, they complement this with a group-wise analysis of femoral shape \textit{trajectories} that included intra-subject, longitudinal data from the OAI. They analyzed this data with the hierarchical model from Sec.~\ref{sec:hierarchical_model} utilizing the integral-based metric discussed in Sec.~\ref{sec:integral-based_metric}.
	
	We compare shape spaces---leave it out? To show that the model is not limited to estimating average, group-level trends, they derived a statistical descriptor for shape trajectories in terms of principal component scores (i.e., the coefficients encoding the trajectories based on principal geodesic modes~\cite{fletcher2004pga}) and used it for trajectory classification.
	
	They determined three groups of shapes trajectories: HH (healthy, i.e., no OA), HD (healthy to diseased, i.e., healthy onset followed by a progression to severe OA), and DD (diseased, i.e., OA at baseline) which where defined by a Kellgren-Lawrence score of grade 0 for all visits, an increase of at least three grades throughout the study, and grade 3 or 4 for all visits, respectively.
	
	For each group, they assembled 22 trajectories, each comprising shapes of all acquired MR images, i.e., at baseline, the 1-, 2-, 3-, 4-, 6-, and 8-year visits.
	
	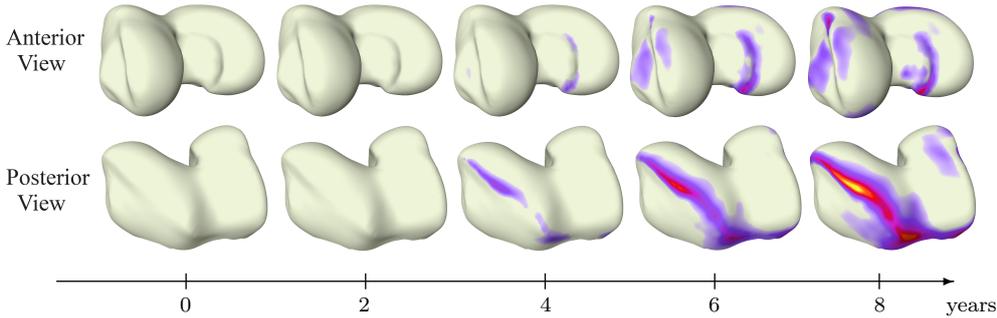
\begin{figure*}[t]
		\begin{center}
			\begin{overpic}[width=1\textwidth,percent,trim={0 0 65mm 0},clip]{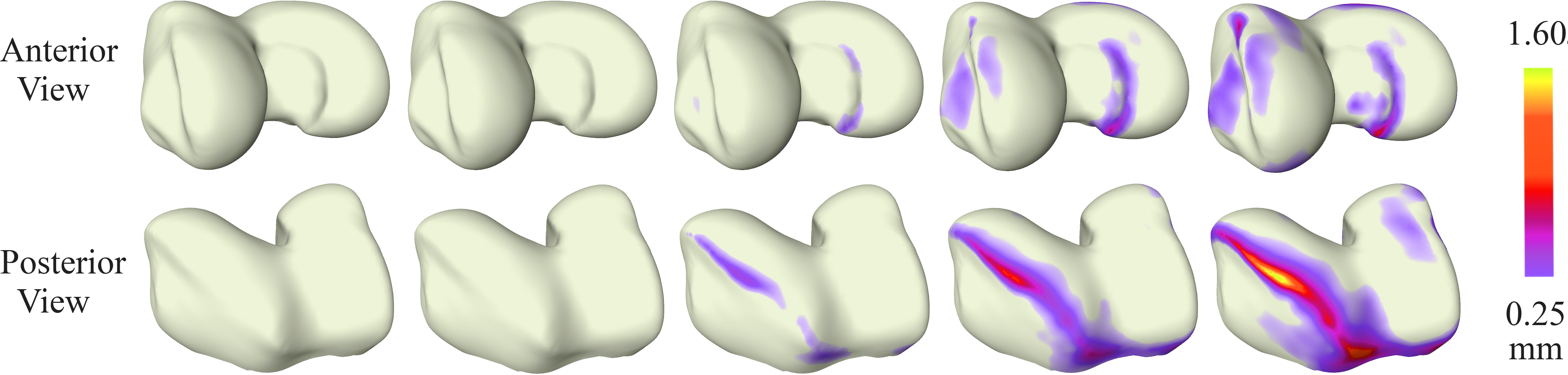}
				\linethickness{0.5pt}
				\put(5,-3){\vector(1,0){0.9\linewidth}}
				\put(87.5,-3){\line(0,1){3pt}\line(0,-1){3pt}}
				\put(87,-6){\footnotesize 8\ \ \ \ \ \ \ \ years}
				\put(18,-3){\line(0,1){3pt}\line(0,-1){3pt}}
				\put(17.4,-6){\footnotesize 0}
				\put(36,-3){\line(0,1){3pt}\line(0,-1){3pt}}
				\put(35.3,-6){\footnotesize 2}
				\put(54,-3){\line(0,1){3pt}\line(0,-1){3pt}}
				\put(53.5,-6){\footnotesize 4}
				\put(71,-3){\line(0,1){3pt}\line(0,-1){3pt}}
				\put(70.4,-6){\footnotesize 6}
			\end{overpic}
			\vspace{2em}
			\caption{Mean of cubic femoral trends of 22 subjects evaluated at five equidistant points. The surface distance to the baseline (value of the computed mean at $t=0$) is color-coded wherever the distance is larger than 0.25 mm according to the color map~(0.25 mm~\ShowColormap~1.6 mm).}
			\vspace{-1em}
			\label{fig:femur_cubic_hierarchical}
		\end{center}
	\end{figure*}
	
	They estimated a hierarchical model for the HD group using cubic B\'ezier curves from $\mathcal{B}^1_3(U)$, $U \subset \Sigma$, to model the individual trends. The degree was chosen because of the results we discussed above.
	
	Time discrete computations were performed based on 2-geodesics under the integral-based metric. Fig.~\ref{fig:femur_cubic_hierarchical} illustrates the estimated group-level trend. Changes along the cartilage plate's ridge are the most prominent OA-related malformations of the femur that are revealed by the determined shape changes. The latter region is known for osteophytic growth. Naturally, the changes are weaker than those of the trajectory in Fig.~\ref{fig:femur_cubic} because the patients that are considered here only ``reached'' KL grade 4 at the end of the study, whereas severe malformations due to several years of grade-4-OA were included in the data that was used in the last section. Similar to the previous experiment, only minute bone remodeling can be observed during the first half of the captured interval, while bone malformations develop more rapidly after four years. Therefore, also this experiment suggests that there are nontrivial phenomena of higher order for which geodesic models are inadequate.
	
	For the classification, the authors first computed the (discrete) mean trajectory of \textit{all} 66 subjects and used it to calculate an approximation $\widetilde{G}$ of the data's Gram matrix; see~\cite{Heeren_ea2018} for details on the latter. For each subject trajectory, they then deduced a 65-dimensional descriptor in the form of the coefficients with respect to the eigenvectors (i.e., PGA modes~\cite{fletcher2004pga}) of $\widetilde{G}$ and trained a simple support vector machine (linear kernel) on the descriptors in a leave-one-out cross-validation setup.
	The percentage of correctly classified trajectories was 64\%. Performing the same experiment with a Euclidean model~\cite{Cootes_ea1995} results in 59\% correct classifications demonstrating the advantage of the Riemannian model over shape spaces that come with the assumption that the data lies in a vector space.

	\subsection{Shape Analysis in Archaeology} \label{sec:archaeology}
	Artifacts---the remnants of tools, weapons, clothes, adornments, and other man-made objects---are the largest and most diversified source of evidence in archaeology, and archaeological analysis has always relied on shape to identify patterns. The ``human dimension'' of artifact creation adds complicated layers of technological, economic, aesthetic, and other social components to shapes. Archaeologists have traditionally utilized \textit{typology} to impose a more or less well-defined order, known as ``typological sequence'' \cite[Ch. 4]{Renfrew2019}, on a series of artifacts. It is hypothesized that this method will isolate the prevalent morphological trends and provide a foundation for extraction and informal analysis of the residual trends that remain unexplained: further geographic detrending through spatial ordering allows for identifying chronological trends; additional temporal ordering may reveal social stratification. Typological ordering (called ``seriation'' in archaeology and ``ordination'' in ecology) has been studied extensively, including computational methods of various complexity (see, e.g., Refs.~\cite{Laxton1989, Madsen1989, Scott1992}). The combinatorial complexity of seriation makes computing the optimal solution impossible even for a small number of artifacts, but archaeological excavations frequently yield tens of thousands. Consequently, no single approach is universally accepted, and manual grouping and ordering of items based on subjective combinations of attributes (artistic/stylistic, technological, utilitarian, etc.) and simplified and abstract graphical representations (drawings) are still prevalent. 
	
	In~\cite{Hanik_ea2023}, the authors proposed regression with Bézier splines in curved shape space for the analysis of parameter-dependent shape trends in artifacts. They also introduced the statistical normalization scheme from Sec.~\ref{sec:normalization} for detrending. As a use case, they investigated the construction principles of ancient Roman and Greek sundials~\cite{Grasshof_ea2016_data}. Therefore, they analyzed the dependence of the shape of the shadow-receiving surface---the ``shadow surface''---of spherical sundials on the latitude of the installation site. For the study, they obtained triangle meshes of the shadow surface of 11 Roman and 3 Greek sundials, all in correspondence; this data is publicly available~\cite{HanikvonTycowicz2022}.
	
	\begin{figure}[ht]
		\centering
		\includegraphics[width=.95\textwidth]{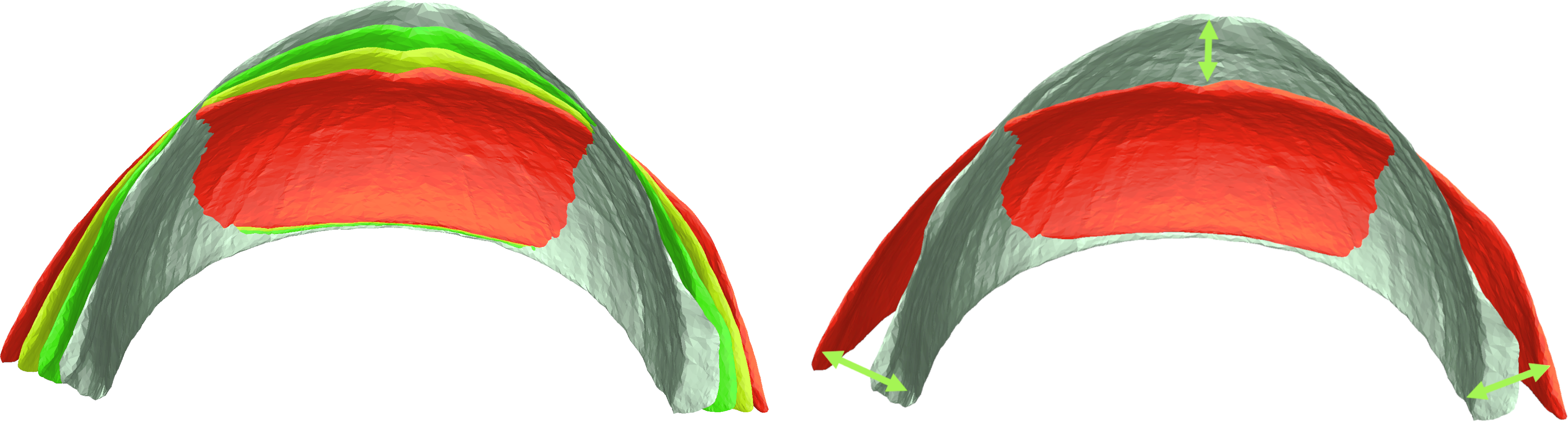}
		\caption{Bending of regressed Roman shadow surfaces, with colors silver, green, yellow, and red ranging (in this order) from the most northern (40.7030\textdegree) to the most southern site (43.3155\textdegree).}
		\label{fig:geodesic_regressions}
	\end{figure}
	
	Using geodesic regression using $\mathcal{B}^1_1(U)$, $U \subset \Sigma$, they find that the shape of the Roman sundials' shadow surface was adapted by means of a latitude-dependent bending; see Fig.~\ref{fig:geodesic_regressions}. They further utilize this knowledge to determine the working latitude of a Roman sundial with an unknown installation site: It can be placed on a latitude only a little north of Rome. They also show how the uncertainty of their prediction can be quantified from the known data and obtain a mean average error of about 80 km in the north-south direction for their prediction.
	
	Finally, the authors used the scheme from Sec.~\ref{sec:normalization} to the shadow surfaces of Greek and Roman at the same latitude. After separately normalizing the shapes at $38,5$\textdegree\ latitude. Since the Greek group only contains 3 samples, the results can only be seen as a first hint. Nevertheless, the 2 mean shapes of the normalized groups clearly look different, which might indicate different construction principles in both regions.
	
	\subsection{An Application in Meteorology: Statistical Analysis of Hurricane Tracks}
	
	Tropical cyclones, also referred to as hurricanes or typhoons, belong to the most supreme natural phenomena with enormous impacts on the environment, economy, and human life. The most common indicator for the intensity of a hurricane is its maximum sustained wind (maxwind), which classifies the storm into categories via the Saffir–Simpson hurricane wind scale. For instance, maxwind $\geq 137$ knots correspond to category 5. The maximal category over a track is called its category. 
	The high variability of tracks and the out-most complexity of hurricanes has led to a huge number of works to classify, rationalize and predict them. We refer to~\cite{REKABDARKOLAEE2019351} for a Bayesian function model,~\cite{Asif2018PHURIEHI} for intensity estimation via machine learning and the overview of recent progress in tropical cyclone intensity forecasting~\cite{RecentProgress}. We remark that many approaches are not intrinsic and use linear approximations. 
	Notable exceptions are the works \citep{SriTrjHur, bauerHomog} that employ an intrinsic Riemannian approach based on the square root velocity framework, albeit only as illustrative examples and without consideration of intensities.
	
\begin{figure}[tb]
		\centering
		\includegraphics[width=.49\textwidth,clip,trim=2cm 11cm 2cm 2cm,cfbox=black .5pt 0pt]{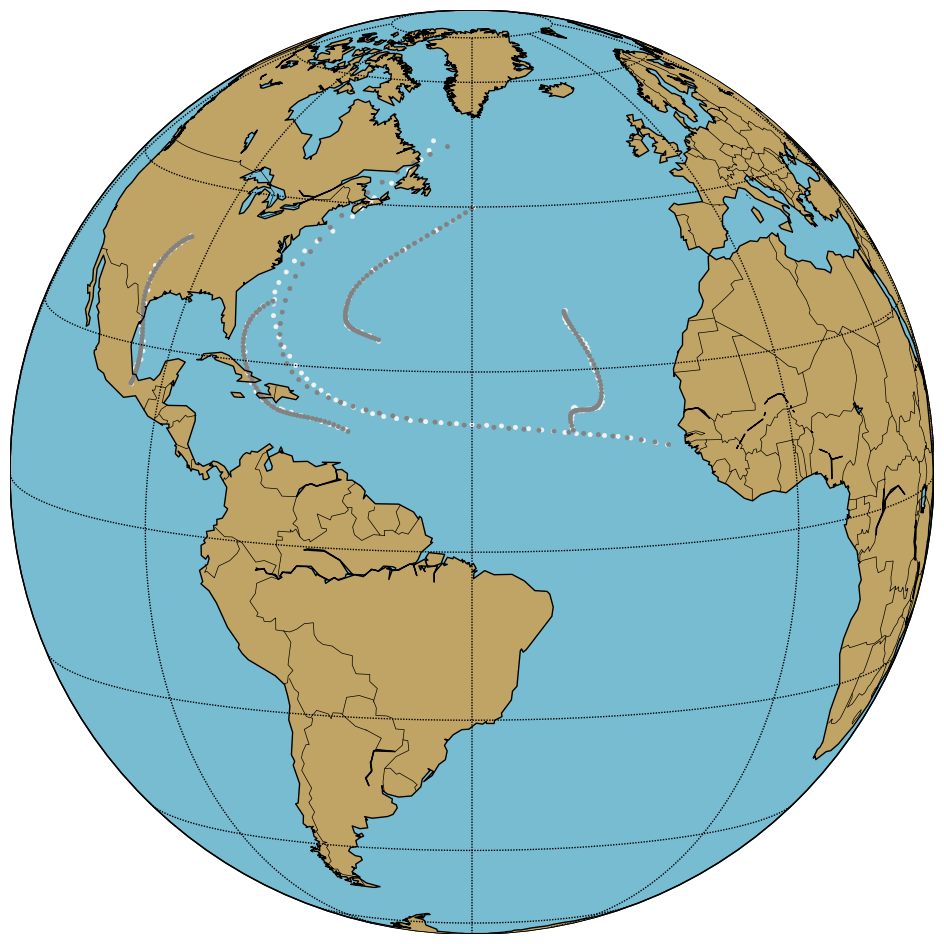}
		\includegraphics[width=.49\textwidth,clip,trim=2cm 11cm 2cm 2cm,cfbox=black .5pt 0pt]{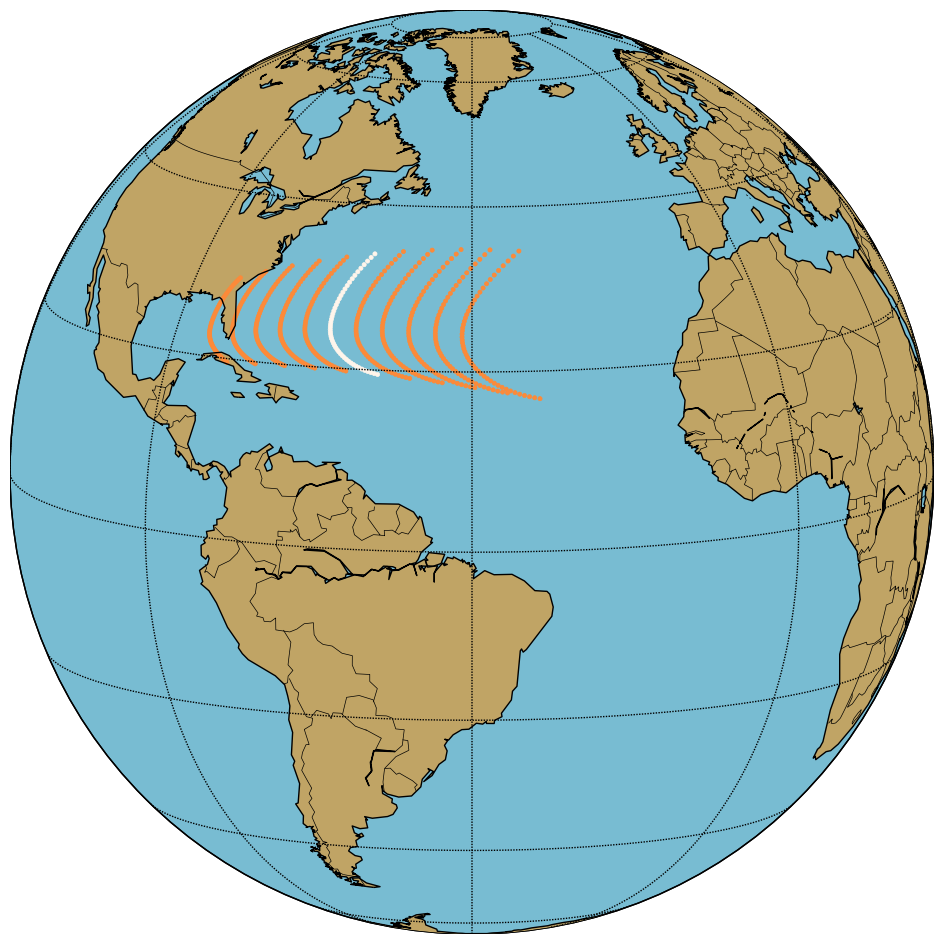}
		\caption{Left: Regressed spline models (grey) of selected Hurricane tracks together with 6 hourly (white) observations. Right: Group-level average track (white) and most dominant mode of variation (orange, sampled within $\pm1.5$ standard deviation).}
		\label{fig:hurricanes}
	\end{figure}

	In~\cite{NavaYazdaniAmbellanHaniketal2023}, the authors analyzed the 2010-2021 Atlantic hurricane tracks (total number 218) recorded in the open-access HURDAT 2 database. The analysis was based on a hierarchical model of the hurricane tracks employing cubic Bézier splines to encode the time-dependent geographical location of each track in $S^2$.
	Empirically 2-segment splines have been found to provide a concise encoding with a very high fidelity in terms of the geometric $R^2$-value.
	Figure~\ref{fig:hurricanes} (left) provides a visualization of selected hurricane tracks showing both the discrete measurements and the corresponding regressed spline representations.
	Subsequently, group-level analysis of the hurricane tracks was carried out via PGA, thus, deriving a summary of the interrelations and characteristic variability in terms of a group-level average track and principal modes of variation. In Fig.~\ref{fig:hurricanes} (right), the average together with the most dominant mode of variation is shown.
	
	Furthermore, they investigated to which extent the (maximal) intensity of a hurricane can be inferred from its trajectory.
	To this end, a support vector machine conditioned on the hierarchical encoding provided by PGA was trained to differentiate between the three intensity classes: (i) tropical storms/depressions (category $<1$), (ii) hurricanes (category $1-3$) with some to devastating damage, and (iii) major hurricanes (category $>3$) with catastrophic damage.
	Using a cross-validation experimental design, a significantly improved discrimination ability of the spline-based representation ($61\%$ and $59\%$ accuracy in average using the Sasakian and integral-based metric, respectively) in comparison to state-of-the-art approaches ($\approx 49\%$ accuracy) that model general, manifold-valued curves was observed.
	This discrepancy in classification performance also prevailed for dimensionality-reduced representations of the general curves, suggesting that the increase in performance can be attributed to the ability of regression schemes to suppress confounding factors such as noise or variances in parameterization.

	\section{Conclusion}
	Ever since its introduction de Casteljau's has proven itself as a fundamental tool in numerous applications ranging well beyond its initially intended use case in digital automotive design. In this overview article, we elaborated on its role in the novel field of geometric data analysis. In particular, we showcased its practicability in obtaining spatiotemporal models that allow for the characterization of parameter-dependent effects in manifold-valued processes. 
	
	The presented geometric statistical tools, i.e., regression and mixed-effects models, are central concepts in statistics and there exist various applications and extensions. While the normalization method we summarized is one example, there are still many fruitful avenues for future developments ranging from geometric model selection to multi-covariate models. Similarly, the presented applications only scratch the surface of potential use cases that will benefit from an intrinsic treatment.   
	
	We further attempted to cover some of the challenges that appear when moving away from the familiar Euclidean domain to manifolds. The lack of properties such as a global system of coordinates and the non-uniqueness of shortest paths often prohibits straightforward generalizations of existing concepts. Also, the lack of closed-form expressions can lead to computational challenges, in particular, for high-dimensional data. 
	
	\section{Acknowledgments}
	Some results in this paper are part one of the author's Ph.D. thesis~\cite{Hanik_phdthesis}.
	M.H.\ is supported by the Deutsche Forschungsgemeinschaft (DFG, German Research Foundation) under Germany´s Excellence Strategy – MATH+: The Berlin Mathematics Research Center, EXC-2046/1 – project ID: 390685689. 
	E.N.\ is supported through the DFG individual funding with project ID 499571814.
	
	\bibliography{reff}
	
\end{document}